\documentclass[a4paper,10pt]{article}
\usepackage{latexsym}
\usepackage{amsmath}
\usepackage{theorem}
\usepackage{epsfig}

\newtheorem{theorem}{Theorem}
\newtheorem{corollary}[theorem]{Corollary}
\newtheorem{definition}{Definition}
\newtheorem{lemma}[theorem]{Lemma}
\newtheorem{proposition}[theorem]{Proposition}

\newcommand{\ZZ}{{\rm\bf Z}}
\newcommand{\RR}{{\rm\bf R}}

\newcommand{\So}{{\mathbf {SO}}(2)}
\newcommand{\EU}{{\rm\bf S}}
\newcommand{\vv}{{\rm\bf v}}
\newcommand{\ww}{{\rm\bf w}}
\newcommand{\dpt}{\displaystyle}
\def\Qed{\hfill\raisebox{.6ex}{\framebox[2.5mm]{}}\medbreak}

\newcounter{contai}
\newenvironment{axiomas}{\begin{list}{{\bf (P\arabic{contai})}}{
\setlength{\leftmargin}{0pt}
\setlength{\labelwidth}{0pt}
\addtolength{\parsep}{10pt}
\usecounter{contai}}}{\end{list}}

\title{Dense heteroclinic tangencies  near a Bykov cycle}

\date{ }

\author{
Isabel S. Labouriau
\quad Alexandre A. P. Rodrigues\\
Centro de Matem\'atica
da Universidade do Porto
\thanks{CMUP (UID/MAT/00144/2013) is funded by FCT (Portugal) with national (MEC) and European structural funds through the programs FEDER, under the partnership agreement PT2020. A.A.P. Rodrigues was supported by the  grants
SFRH/BD/28936/2006 and SFRH/BPD/84709/2012 of FCT.}\\
 and
Faculdade de
Ci\^encias, Universidade do Porto \\
Rua do Campo Alegre,
687, 4169-007 Porto, Portugal \\
 islabour@fc.up.pt \quad alexandre.rodrigues@fc.up.pt }
\begin{document}

\maketitle


\bigbreak
\textbf{2010 --- AMS Subject Classifications}

  {Primary: 34C28

   Secondary: 34C37, 37C29, 37D05, 37G35}

\bigbreak

\maketitle

\begin{abstract}
This article presents a mechanism  for the coexistence of hyperbolic and non-hyperbolic dynamics arising in a neighbourhood of a Bykov cycle where trajectories turn in opposite directions near the two nodes --- we say that the nodes have different chirality.
We show that in  the set of vector fields defined on a  three-dimensional manifold, there is a class where tangencies of the invariant manifolds of two hyperbolic saddle-foci occur densely.
The class is defined by the presence of the Bykov cycle, and by a condition on the parameters that determine the linear part of the vector field at the equilibria.
This has important consequences: the  global dynamics is persistently dominated by heteroclinic tangencies and by Newhouse phenomena, coexisting with hyperbolic dynamics arising from  transversality.
The coexistence gives rise to linked suspensions of Cantor sets,  with hyperbolic and non-hyperbolic dynamics, in contrast with the case where the nodes have the same chirality.

We illustrate our theory with an explicit example where tangencies arise in the unfolding of a symmetric vector field on the three-dimensional sphere.
\end{abstract}

\section{Introduction}
\begin{figure}
\begin{center}
\includegraphics[height=4cm]{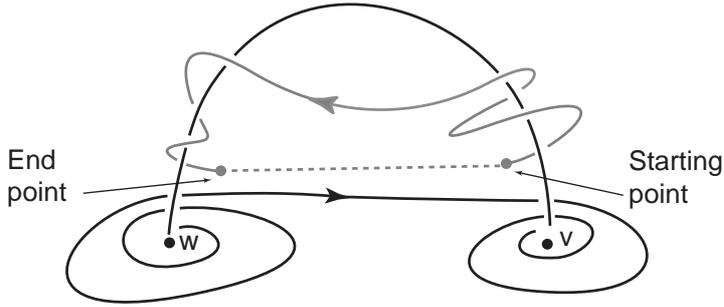}
\end{center}
\caption{\small Bykov cycle with saddle-foci of different chirality. The starting point of a nearby  trajectory is joined to its end point  forming a loop. Arbitrarily close to the cycle there are trajectories whose loop is not linked to the cycle. This happens because near each  saddle-focus
trajectories turn around the one-dimensional connection  in   opposite directions.}
\label{Bykov2}
\end{figure}

Consider a differential equation in  a three-dimensional manifold having a heteroclinic cycle that consists of two saddle-foci of different Morse indices whose one-dimensional invariant manifolds coincide and whose two-dimensional invariant manifolds intersect transversely.
There  are two different possibilities for the geometry of the flow around the cycle, depending on the direction trajectories turn around the  heteroclinic connection of one-dimensional  invariant manifolds. 
The two cases give rise to different dynamics, but the distinction is usually not made explicitly in the literature. 
This article is concerned with the case when, 
near the two saddle-foci, trajectories wind in opposite directions around the heteroclinic connection of one-dimensional invariant manifolds  --- the two nodes have \emph{different chirality} as in Figure \ref{Bykov2}.

The dynamics around this type of cycle was first studied by V.V. Bykov \cite{Bykov99,Bykov}, with the implicit assumption of different chirality. 
He has obtained an open class containing a dense set of vector fields exhibiting tangencies of the two-dimensional invariant manifolds --- see also \cite[Th. 5.33]{HS}. 
Bykov also described bifurcations occurring when the  structurally unstable one-dimensional connection is broken.

In the present article we highlight that the orientation of the flow around the  structurally unstable connection has profound effects on the dynamics near the cycle.
We refine and clarify the key ideas of  the analysis of the unperturbed system of \cite{Bykov99, Bykov}  and we  explore properties of the maximal invariant set that emerges near the Bykov cycle. 
After recovering Bykov's one-pulse heteroclinic tangencies, we show the existence of multi-pulse connections, occurring along trajectories that follow the original cycle an arbitrary number of times.
We also show that the non-hyperbolic set containing heteroclinic tangencies coexists with the suspension of uniformly hyperbolic horseshoes.
Although each individual tangency may be eliminated by a small perturbation, another tangency is created nearby, while the hyperbolic set persists.
By construction, the single-round periodic solutions and the single-round heteroclinic trajectories found in \cite[Th. 3.1 and 3.2]{Bykov} lie inside the suspended horseshoes found here.
An explicit vector field where the nodes have different chirality is also constructed here.

Tangencies of invariant manifolds are associated to
Newhouse phenomena:
bifurcations leading
to the birth of infinitely many  asymptotically stable periodic solutions \cite{Newhouse2, Newhouse1}.
Such tangencies have been recognised as a mechanism for instability and lack of hyperbolicity
in surface diffeomorphisms.
 One of the first results in this direction
was established by Gavrilov and Shilnikov \cite{GS} for two-dimensional maps. To understand this phenomenon it is important to study  the variety of dynamical behaviour associated to the creation and destruction of tangencies.
Although periodic attractors arising from tangencies
have quite large periods and small basins of attraction,
an infinite number of them
may change the character of the chaotic dynamics.
Around a Bykov cycle with nodes of different chirality, 
suspended uniformly hyperbolic horseshoes  coexist  with tangencies and are not separated as a whole, in sharp contrast to what is expected of attractors that are either hyperbolic or Lorenz-like.
\bigbreak

Our work also forms part of a program started by Glendinning and Sparrow \cite{T-points Glendinning}  in the eighties and by Bykov \cite{Bykov99, Bykov} in the late nineties,  addressing the
systematic study of the dynamics near networks of equilibria  whose linearisation has a conjugate pair of  non real eigenvalues.  
Recently, there has been a renewal of interest in this type of heteroclinic bifurcation in the reversible and conservative cases \cite{KLW, LR, LTW}. 
 They have been discussed in the context of reversible divergence-free systems under the name T-points in \cite{KLW, LTW}.  See also the survey \cite{HS}. 
 These works describe the types of nonwandering dynamics and bifurcations nearby. The discussion about the orientation of the flow around the one-dimensional connection has not yet appeared in the literature, but has been implicitly  considered by several authors.

Without the condition of different chirality the density of heteroclinic tangencies near the cycle does not hold, as we proceed to discuss.
In the cases  \cite{ACL NONLINEARITY, ALR, KLW, LR, LR2, LR3}, the authors assumed, sometimes implicitly, that  trajectories wind in the same direction around the one-dimensional connection in the neighbourhood of the saddle-foci.

An example of a heteroclinic network involving several hyperbolic equilibria,
where trajectories switch around the different cycles of the network, is described by Aguiar \emph{et al} \cite{ACL NONLINEARITY}.
A symmetry reduction argument yields a quotient network with
two saddle-foci of different Morse indices reminiscent of those studied by Bykov \cite{Bykov}.
Under the assumption that near the two equilibria trajectories wind in the same direction around the one-dimensional heteroclinic connection,
it is shown that each small tubular neighbourhood of the network contains suspended horseshoes  \cite{ACL NONLINEARITY, LR}.
In this case, some of the connections in the network arise from transverse intersections of stable and unstable manifolds of equilibria.
Rodrigues \cite{Rodrigues3}
has proved that Lebesgue almost all solutions do not remain near the cycle for all time, although they may return to the cycle after an excursion away from it.

The recent work by Knobloch \emph{et al} \cite[\S 6.2.1]{KLW} confirms, using Lin's method, that heteroclinic tangencies near the cycle are rare.  
All  the analysis has been done assuming the same chirality, which comes for free in the reversible context induced by the Michelson system.
Instead of restricting the flow to
the wall of the cylinder (as we do in Theorem \ref{prop12}), their conclusions, like Bykov's,  are supported by the study of spirals in another cross section. They proved that the spirals corresponding to the invariant manifolds have at most two tangential intersections in total. 
Under the hypothesis of same chirality, the non-wandering dynamics near the Bykov cycle is dominated by hyperbolic horseshoes, conjugate to a full shift over a finite alphabet, that accumulate on the cycle. See \cite{LR, Rodrigues3}.

Heteroclinic tangencies arise near cycles with nodes of the same chirality in systems close to symmetry
\cite{LR2015}, but in the latter case the trajectories always make an excursion far from the cycle.
In contrast, tangencies found in this article concern trajectories that remain close to the cycle for all time.

Needless to say that, for an experimentalist, the analysis of these phenomena is of great importance. For instance, our results can be useful to locate cocoon bifurcations near a Bykov cycle, a phenomenon observed by Lau in the Michelson system --- details in \cite{DumortierEtAl}.

\bigbreak

In summary, we establish here that chirality is a topological invariant for cycles in networks that determines whether or not the tangencies accumulate on the cycle.  Many questions remain for future work. One obvious one is to construct a Bykov cycle with two saddles with different divergence sign. In this case, it might be possible to show the existence of mixed dynamics  \cite{Delshams} in the neighbourhood of the Bykov cycle without additional assumptions on the reversibility. Questions concerning the global bifurcations  associated to symmetry breaking from regular dynamics to chaos would also be of great interest.

\bigbreak
\textbf{Structure of the article:}
We study the set of non-wandering points near a Bykov cycle in $\textbf{S}^3$ and show the existence of heteroclinic tangencies.
In Section~\ref{SecDescription}, after recalling some preliminary definitions, we state the main results, Theorems~\ref{teoremaReverteRotacao} and \ref{TeoremaTransversal}.
Their dynamical consequences are discussed in Section~\ref{NewDiscussion}.
We establish the notation and framework for the proof of Theorem \ref{teoremaReverteRotacao} in Section \ref{SectionLocalDynamics} where we linearise the vector field around each equilibrium obtaining a geometrical description of the way the flow transforms a curve of initial conditions lying across the stable manifold of an equilibrium.
Section~\ref{SectionT-pointHetero} contains   more precise statements of
Theorems~\ref{teoremaReverteRotacao} and~\ref{TeoremaTransversal}
and their proofs.
 This is followed  in Section~\ref{numerics} by an example of a symmetric vector field with these properties, illustrated by numerical simulations.We finish the article with a short discussion of the implicit assumption of different chirality in Bykov's articles \cite{Bykov99,Bykov}.

\section{Description of the problem}\label{SecDescription}
\subsection{Preliminaries}

Let $f$ be a smooth vector field defined on a 3-dimensional smooth manifold. Given two equilibria $p$ and $q$  of $\dot{x}=f(x)$, a \emph{heteroclinic connection} from $p$ to $q$, denoted  $[p\rightarrow q]$, is a
flow-invariant subset of $W^u(p)\cap W^s(q)$.
In this article we consider mostly 1-dimensional, sometimes 2-dimensional, connections between equilibria.
We are interested in \emph{heteroclinic cycles} associated to two hyperbolic equilibria $p$ and $q$: the  set consisting of the equilibria and two heteroclinic connections $[p \rightarrow q] , [q \rightarrow p]$. Sometimes we refer to the equilibria on the  cycle as \emph{nodes}.
More general information may be found in  \cite{AC, Field book}.

The dimension of the  unstable manifold of a hyperbolic equilibrium is called the \emph{Morse index} of the equilibrium.
A saddle-focus $p$ is an equilibrium of  $\dot{x}=f(x)$  where the spectrum of $df(p)$ has two complex non-real eigenvalues $\alpha \pm i\omega$ and one real eigenvalue $\beta$ with $\alpha\beta<0$.
A  \emph{Bykov cycle} is a heteroclinic cycle associated to two hyperbolic saddle-foci with different Morse indices, in which the one-dimensional invariant manifolds coincide and the two-dimensional invariant manifolds have a transverse intersection.
The presence of the two connections in a Bykov cycle implies the existence of infinitely many subsidiary connections following the original cycle -- see Labouriau and Rodrigues \cite{LR}.

\subsection{Hypotheses}
\label{description}
\begin{figure}
\begin{center}
\includegraphics[height=5cm]{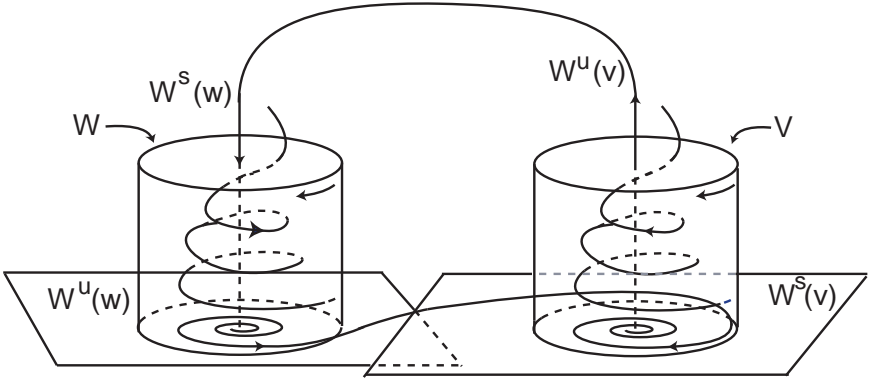}
\end{center}
\caption{\small Geometry near a Bykov cycle with saddle-foci of different chirality.
A trajectory starting in a neighbourhood  $V$ of $\vv$ turns around  the connection  $[\vv \rightarrow \ww]$ in one direction. After the trajectory arrives at a neighbourhood $W$ of $\ww$ it  turns around    $[\vv \rightarrow \ww]$ in the opposite  direction.}
\label{orientations}
\end{figure}

Our object of study is the dynamics around a special type of Bykov cycle, for which we give a rigorous description here. Specifically, we study a $C^2$--vector field $f$ on a manifold diffeomorphic to the three-sphere $\EU^3 =\{X=(x_1,x_2,x_3,x_4) \in \RR^4: ||X||=1\}$ whose flow has the following properties (see Figure \ref{orientations}):

\begin{axiomas}
\item\label{H1} There are two hyperbolic saddle-foci $\vv$ and $\ww$. The eigenvalues of $df_X$ are:
\begin{itemize}
\item[(a)] $- C_\vv \pm \alpha_\vv i$ and $E_\vv$  where $C_\vv$, $E_\vv$ and $\alpha_\vv$ are positive, for $X= \vv$;
\item[(b)] $ E_\ww \pm \alpha_\ww i$ and $-C_\ww$ where $C_\ww$, $E_\ww$ and $\alpha_\ww$ are positive, for $X= \ww$.
\end{itemize}

\item \label{H2} There is a heteroclinic cycle $\Gamma$ consisting of $\vv$, $\ww$ and two one-dimensional heteroclinic connections $[\vv \rightarrow \ww]$ and $[\ww \rightarrow \vv]$.

\item \label{H3} At the heteroclinic connection $[\ww \rightarrow \vv]$, the two-dimensional manifolds $W^u(\ww)$ and $W^s(\vv)$ meet transversely.

There are two different possibilities for the geometry of the flow around $\Gamma$,
depending on the direction trajectories turn around the heteroclinic connection $[\vv\to\ww]$.
To make this rigorous, we need some new concepts. Let $V$
and $W$ be small disjoint neighbourhoods of $\vv$ and $\ww$ with disjoint boundaries $\partial V$ and $\partial W$, respectively.
Trajectories starting at $\partial V$ near $W^s(\vv)$ go into the interior of $V$ in positive time, then follow the connection $[\vv\rightarrow\ww]$, go inside $W$, and then come out at $\partial W$
as in Figure~\ref{orientations}.
Let $\varphi$ be a piece of trajectory like this from $\partial V$ to $\partial W$.
Now join its starting point to its end point by a line segment as in Figure~\ref{Bykov2}, forming a closed curve, that we call the  \emph{loop} of $\varphi$.
For generic starting points, the loop of $\varphi$ and the cycle $\Gamma$ are disjoint closed curves. We say that the two saddle-foci $\vv$ and $\ww$ in $\Gamma$ have the same \emph{chirality} if the loop of every trajectory is linked to $\Gamma$ in the sense that the two closed curves cannot be disconnected by an isotopy. Otherwise, we say that $\vv$ and $\ww$ have different chirality: given the neighbourhoods $V$ and $W$, it is always possible to find a trajectory $\varphi$ going through $V$ and $W$ whose loop is not linked to $\Gamma$. Then our assumption is:

\item \label{H4} The saddle-foci $\vv$ and $\ww$ have different chirality.
\end{axiomas}

The transverse intersection of $W^u(\ww)$ and $W^s(\vv)$ of (P\ref{H3}) persists under $C^1$--perturbations, whereas the connection $[\vv \rightarrow \ww]$ does not, unless there is some special property, like symmetry.

Property (P\ref{H4}) is persistent under isotopies: if it holds for $f$, then it is still valid in continuous (not necessarily smooth) one-parameter families containing it, as long as there is still a connection.  This is particularly important  when we consider unfoldings of $f$.

\subsection{Statement of the main results}

The main result in this article guarantees that heteroclinic tangencies {are ubiquitous in systems having} Bykov cycles with different chirality.
We prove that these tangencies occur  {densely,  when the parameters that determine the linear part of the vector field at the equilibria lie in a set of   full Lebesgue measure. The tangencies} 
%
 lie near the cycle, in contrast to the findings of \cite{KLW,  LR2, LTW, LR3} for cycles with the same chirality, where either they  do not appear at all or, as in \cite{LR2015}, they appear far from the cycle.

Given two disjoint neighbourhoods $V$ of $\vv$
and $W$  of $\ww$ with disjoint boundaries, $\partial V$ and $\partial W$, respectively,
consider a point in $[\vv\rightarrow\ww]\cap \partial V$ and a  neighbourhood in $\partial V$ of this point  that is also a cross-section to $f$.
Saturating the cross-section by the flow, one obtains a flow-invariant tube joining $V$ to $W$   { that contains} the connection in its interior.
A similar flow-invariant tube may be obtained around the connection $[\ww\rightarrow\vv]$ joining $W$ to $V$.
We call the union of $V$ and $W$ with these tubes,  a \emph{tubular neighbourhood}
of the Bykov cycle. The next results hold for vector fields in $\EU^3$ under  the $C^k$ topology, for $k\ge 2$.

\begin{theorem}
\label{teoremaReverteRotacao}
There is an open set $\mathcal{C}$ of vector fields satisfying (P\ref{H1})--(P\ref{H4}) such that for any $f\in\mathcal{C}$, and any tubular neighbourhood $U$ of the Bykov cycle,
 there are vector fields on $\EU^3$  arbitrarily close to $f$,  with a Bykov cycle in $U$ with the same properties, for which $W^u(\ww)$ and $W^s(\vv)$ have a tangency inside $U$.
\end{theorem}

A more precise statement of this result will be given in Section \ref{SectionT-pointHetero}. The open set $\mathcal{C}$ corresponds to an open condition on the eigenvalues of the linearisation of the vector field that, by  Thom's Transversality Theorem, defines an open set in the $C^2$ topology.
The set of vector fields with a heteroclinic tangency is  dense in $\mathcal{C}$.

\begin{definition}
Let $W$ be a small neighbourhood of $\ww$ and let $\Sigma\subset W$ be a cross-section to the flow meeting $W^u(\ww)$. A one-dimensional connection
$[\ww \rightarrow \vv]$ that meets $\Sigma$ at precisely $n\in{\rm\bf N}$ points is called a $n$-pulse heteroclinic connection $[\ww \rightarrow \vv]$ with respect to $\Sigma$.
\end{definition}

The tangencies of invariant manifolds coexist with transverse intersections, giving rise to a hyperbolic structure similar to that
obtained when (P\ref{H4}) does not hold as in \cite{ACL NONLINEARITY,LR}.

\begin{theorem}\label{TeoremaTransversal}
There is an open set  $\mathcal{E}$ of vector fields in $\EU^3$ satisfying (P\ref{H1})--(P\ref{H4}) for which
any tubular neighbourhood $U$ of the Bykov cycle $\Gamma$ contains the following:
\begin{enumerate}
\item \label{item6}
 trajectories in $U\backslash \Gamma$ that remain in $U$ for all time;
\item \label{item4}
at least one $n$-pulse heteroclinic connection $[\ww\to\vv]$ for each $n \in \textbf{N}$;
\item \label{item5}
a  cross-section $S\subset U$ containing a set of points such that at these points the dynamics of the first return to  $S$ is  uniformly hyperbolic and conjugate to a full shift over a finite number of symbols.
 This set accumulates on the cycle.
\end{enumerate}
Moreover,  the set  $\mathcal{E}$ meets the set $\mathcal{C}$ of
Theorem~\ref{teoremaReverteRotacao}  on a non-empty open set of vector fields.
\end{theorem}

\begin{corollary}\label{corolarioCoexiste}
For vector fields satisfying (P\ref{H1})--(P\ref{H4}) and in a dense subset of the open set
$\mathcal{E}\cap\mathcal{C}$
tangencies of the two-dimensional invariant manifolds of the saddle-foci
coexist with transverse intersections.
\end{corollary}

\section{Dynamical Consequences}\label{NewDiscussion}

We will show that transverse and tangent heteroclinic connections of two-dimensional invariant manifolds coexist  near a Bykov cycle with  nodes of different chirality.
In this section we explore the consequences of this result.

Our analysis has been restricted to $\EU^3$, the lowest possible phase space dimension
in which Bykov cycles may occur.
Extension to  higher dimensions may be possible using
the heteroclinic centre manifold theorem, Lin's method and similar techniques \cite{KLW, LTW}; in this article we have not attempted to do so.

\subsection{Chirality}
The conclusions of Theorem~\ref{TeoremaTransversal} do not depend on chirality and are the same as  those for Bykov cycles of the same chirality in \cite{ACL NONLINEARITY, LR,LR2, LTW}: any tubular neighbourhood of the cycle contains
trajectories that remain on it for all time forming infinitely many suspended horseshoes.

The main difference due to property (P\ref{H4})  is that any neighbourhood of a Bykov cycle with different chirality contains nontrivial and irremovable subsets with hyperbolic dynamics but it is not exhausted by them. The non-hyperbolicity takes place everywhere.

\begin{corollary}
There is an open set $\mathcal{C}$ of vector fields satisfying (P\ref{H1})--(P\ref{H4}) such that arbitrarily close to any $f\in\mathcal{C}$, there is a vector field on $\EU^3$ with a Bykov cycle with the same properties,
where the non-uniformly hyperbolic dynamics cannot be separated by an isotopy
 from the maximal hyperbolic set that appears  in any tubular neighbourhood of the cycle.
 \end{corollary}

Equations with symmetry having cycles whose nodes have the same chirality may also exhibit tangencies of invariant manifolds, but they arise  far from the Bykov cycle, see  \cite{LR2, LR2015}.

\subsection{Tangency}
In the context of  dissipative diffeomorphisms containing homoclinic points, Newhouse \cite{Newhouse2, Newhouse1} introduced the term \emph{wild attractor} for non-uniformly hyperbolic sets whose invariant manifolds have a tangency. He reported what happens in a one-parameter
unfolding, when a homoclinic tangency splits, and discovered nontrivial, transitive and hyperbolic sets whose stable and unstable invariant manifolds have an irremovable nondegenerate tangency --  although one tangency may be removed by a small perturbation of the system, one cannot avoid the appearance of new ones.
Fat hyperbolic sets, in the sense of Bowen \cite{Bowen, Bowen75}, exist for diffeomorphisms $C^2$--close to any diffeomorphism with a homoclinic tangency --- see Palis and Takens \cite[\S 4]{PT}. The open regions where diffeomorphisms with homoclinic tangencies are dense are called \emph{Newhouse regions}.
Newhouse's results on homoclinic tangencies may be adapted to the case studied here, ensuring
the multiplicity of nearby stable solutions.

\begin{corollary}\label{corollaryPeriodDoubling}
Let $f_\lambda$ be a one-parameter family of vector fields in the open set $\mathcal{C}$ of Theorem~\ref{teoremaReverteRotacao} and suppose in addition that $C_\vv>E_\vv$ and $C_\ww>E_\ww$.
Then there are period doubling sequences for the first return map to any  cross section to the Bykov cycle.
Moreover, there are persistent heteroclinic tangencies of the invariant manifolds of periodic solutions and infinitely many periodic sinks nearby.
\end{corollary}

\textbf{Proof:}
This follows from  results  by Mora and Viana \cite{MV}, Palis and Takens \cite{PT} and Yorke and Alligood \cite{YA}, to which the reader is also referred for more details on the bifurcation sequences giving rise to these dynamical properties.
\Qed

Near a heteroclinic tangency there is no dominated  splitting of the tangent space into  stable and  unstable subspaces.
In the $C^1$--topology, for two-dimensional maps, the existence of tangencies and
the absence of dominated splittings are  {synonyms} \cite{Wen}.

 In the conservative and reversible settings, it is known that there are coexisting Newhouse regions in
which a dense set of maps possess simultaneously infinitely many asymptotically stable, saddle, completely unstable and elliptic periodic orbits --- see \cite{Delshams, LS}  {and references therein.}

The description of the set of all solutions that lie near a non-transverse intersection of invariant manifold becomes complicated.
 As reported in Gonchenko \emph{et al} \cite{GST}, the main
source of the difficulty is that arbitrarily small perturbations of any
system with the simplest tangency may lead
to the creation of new tangencies of arbitrarily high
orders, and to the birth of periodic trajectories of arbitrarily high orders of degeneracy.

\subsection{The spiralling set}

A transverse intersection of  two-dimensional invariant manifolds $W^u(\ww)$ and $W^s(\vv)$ cannot be removed by a small smooth perturbation.
Near it the set of all solutions that never leave a neighbourhood of the heteroclinic cycle
 has a non-trivial structure as in  Shilnikov \cite{Shilnikov_67}: it defines a locally maximal non-trivial hyperbolic set, and admits
a complete description in terms of symbolic dynamics.
When the non transverse one-dimensional connection is broken, for instance by forced symmetry-breaking as in \cite{LR, LR3}, finitely many of these horseshoes persist  but their number may be arbitrarily large.

The existence of hyperbolic horseshoes near the cycle gives some information about the topological entropy of the first return map \cite{Katok2}. Although, near the horseshoes, there are irremovable sinks  with zero topological entropy, we may still conclude that:

\begin{corollary}
For $f$ in the open set $\mathcal{C}$ of Theorem~\ref{teoremaReverteRotacao},
the topological entropy associated to the first return map restricted to a cross section $\Sigma$ to the Bykov cycle  is positive.
\end{corollary}

The previous results  {imply} that there are trajectories with one positive Lyapunov exponent.
The global attractor  also contains infinitely many sinks with long
periods and narrow basins of attraction, arising from the heteroclinic tangencies described in Theorem~\ref{teoremaReverteRotacao}.
The maximal transitive set surrounds attracting periodic solutions that accumulate on the original cycle and coexist with hyperbolic horseshoes.
The properties of the attractor are similar to the quasi-stochastic attractors studied by Gonchenko \emph{et al} \cite{GSTchaos}.

Homoclinic tangencies of arbitrary orders are considered by Gonchenko and Li \cite{Gonchenko} together with sufficient conditions for the existence of nontrivial hyperbolic sets containing infinitely many hyperbolic horseshoes.
Their proof uses symbolic dynamics and refers to  arguments of Katok \cite{Katok2}.

\section{Local Dynamics near each saddle-focus}
\label{SectionLocalDynamics}
We follow the standard procedure for describing the dynamics near a heteroclinic network. We construct return maps
defined on various cross-sections and analyse the dynamics
 by composing them in an appropriate order to obtain Poincar\'e maps modelling the dynamics near the Bykov cycle.

In this section, we establish local coordinates near the saddle-foci $\vv$ and $\ww$ and define some notation that will be used in the rest of the paper. The starting point is an application of Samovol's Theorem \cite{Samovol} to linearise the flow around the equilibria and to introduce cylindrical coordinates around each saddle-focus.
These are used to define neighbourhoods whose boundaries are
transverse to the linearised flow. For each saddle, we obtain the expression of the local map that sends points in the boundary where the flow goes in, into points in the boundary where the flows goes out. Finally, we establish a convention for the transition maps from one neighbourhood to the other. When we refer to the stable/unstable manifold of an equilibrium point, we mean the \emph{local} stable/unstable manifold of that equilibrium.

\subsection{Linearisation near the saddle-foci}
\label{linearization}
By Samovol's Theorem \cite{Samovol}, around the saddle-foci, the vector field $f$ is $C^1$--conjugate to its linear part, since there are no resonances of order 1.
In cylindrical coordinates $(\rho ,\theta ,z)$ the linearisations at $\vv $ and $\ww $ take the form, respectively:
\begin{equation}\label{localMapv}
\left\{
\begin{array}{l}
\dot{\rho}=-C_{\vv }\rho \\
\dot{\theta}=\alpha_{\vv } \\
\dot{z}=E_{\vv }z
\end{array}
\right.
\qquad
\left\{
\begin{array}{l}
\dot{\rho}=E_{\ww }\rho \\
\dot{\theta}= - \alpha_{\ww } \\
\dot{z}=-C_{\ww }z .
\end{array}
\right.
\end{equation}

\begin{figure}
\begin{center}
\includegraphics[width=12cm]{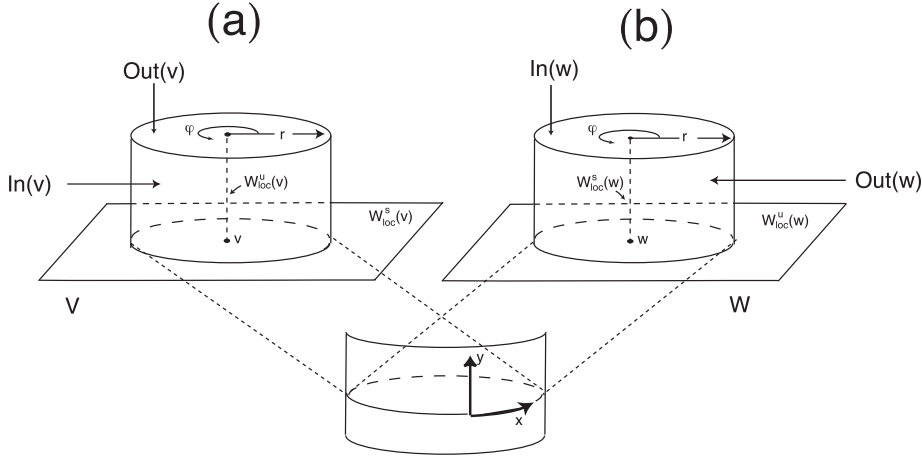}
\end{center}
\caption{\small Parametrisation of the cylindrical neighbourhoods of the  saddle-foci.
(a)
The flow goes into the cylinder $V$ transversely across the wall $In(\vv) \backslash W^s(\vv)$
and leaves it transversely across its top $Out(\vv)$.
 (b)
 The flow goes into the cylinder $W$ transversely across its top  $In(\ww) \backslash W^s(\ww)$ and leaves it transversely across the wall $Out(\ww)$.  Inside both cylinders, the vector field is linear.
 Cylinder tops are parametrised in polar coordinates $(r,\varphi)$, cylinder walls in coordinates
 $(x,y)$, with angular coordinate $x$.}
\label{FigvwSections}
\end{figure}

We consider
cylindrical neighbourhoods of $\vv $ and $\ww $ in $\textbf{S}^3$ of radius $\varepsilon>0$ and height $2\varepsilon$ that we denote by $V$ and $W$, respectively. Their boundaries consist of three components (see Figure \ref{FigvwSections}):
\begin{itemize}
\item
The cylinder wall parametrised by
$x\in \RR\pmod{2\pi}$ and $|y|\leq \varepsilon$ with the usual cover:
$$(x,y)\mapsto (\varepsilon ,x,y)=(\rho ,\theta ,z).$$ Here $y$ represents the  height of the cylinder and $x$ is the angular coordinate, measured from the point $x=0$ in the connection $[\ww\to\vv]$.
\item
Two disks, the top and the bottom of the cylinder.
We assume the connection $[\vv\to \ww]$ goes from the top of one cylinder to the top of the other, and we take a polar covering of the top disk:
$$(r,\varphi )\mapsto (r,\varphi ,\varepsilon )=(\rho ,\theta ,z)$$
where
$0\leq r\leq \varepsilon $ and $\varphi \in \RR\pmod{2\pi}$.
\end{itemize}
On these cross sections, we define the return maps to study the dynamics near the cycle.

\subsection{Coordinates near $\vv$ and $\ww$}
Consider the cylinder wall  near $\vv$, that  {locally} meets  $W^{s}(\vv )$ on the circle parametrised by $y=0$.
The top part  $y\ge 0$ of the cylinder wall  near $\vv$ is denoted by $In(\vv)$. Trajectories starting at interior points of $In(\vv)$ go into the cylinder in positive time and come out at the cylinder top,  denoted $Out(\vv)$.
Trajectories starting at interior points of $Out(\vv)$ go inside the
cylinder in negative time. After linearisation,
in these coordinates, the manifold $W^{u}(\vv )$ is the $z$--axis, intersecting $Out(\vv)$ at the
origin of coordinates.

Reversing the time, we get dual results for $\ww$. After linearisation, $W^{s}(\ww )$ is the $z$--axis, intersecting the top, $In(\ww)$, of the cylinder at the origin of its coordinates.
Trajectories starting at  interior points of $In(\ww)$
go into $W$ in positive time.

Trajectories starting at interior points of the cylinder wall  $Out(\ww)$
go into $W$ in negative time.
The set $Out(\ww)\cap W^{u}(\ww )$ is parametrised by $y=0$.
Trajectories that start at
$In(\ww)\backslash W^s(\ww )$  leave the cylindrical
neighbourhood $W$ at $Out(\ww)$ .

\subsection{Transition Maps}
In the rest of this paper, we study the Poincar\'e first return map on the boundaries defined above.
Consider the  transition maps
$$
\Psi_{\vv,\ww}:Out(\vv) \longrightarrow In(\ww)
\qquad
\text{and} \quad \Psi_{\ww,\vv}:Out(\ww) \longrightarrow In(\vv).
$$
The map $\Psi_{\ww,\vv}$ can be seen as a rotation by an angle $\alpha$. We use $\alpha \equiv \frac{\pi}{2}$, that simplifies the expressions used.

As in Bykov \cite{Bykov} and Homburg and Sandstede \cite{HS}, after a rotation and a uniform rescaling of the coordinates, we may assume without loss of generality that $\Psi_{\vv, \ww}$ is the linear map
$\Psi_{\vv, \ww}(x,y)=\left(a x,y/a \right)$ with $a\ge 1$.

\subsection{Local maps near $\vv$ and $\ww$}
The flow is transverse to the above cross sections and moreover the boundaries of $V$ and of $W$ may be written as the closures of the disjoint unions $In(\vv) \cup Out (\vv)$ and  $In(\ww) \cup Out (\ww)$, respectively. The trajectory of  the point $(x,y)$ in $In(\vv) \backslash W^s(\vv)$ leaves $V$ at $Out(\vv)$ at:
$$
\begin{array}{c}
\Phi_{\vv }(x,y)=(c_{1}y^{\delta_\vv},-g_\vv \ln y+x+c_{2})=(r,\varphi)
\end{array},
$$
where $\delta_\vv$ is the \emph{saddle index} of $\vv$,
$$
\delta_\vv=\frac{C_{\vv }}{E_{\vv}}>0,
\quad
c_{1}=\varepsilon ^{1-\delta_\vv}>0,
\quad
g_\vv=\frac{\alpha_{\vv }}{E_{\vv}}>0
\quad
\text{ and }
\quad
c_{2}=g_\vv\ln (\varepsilon ).
$$

Similarly, points $(r,\varphi)$ in $In(\ww) \backslash W^s(\ww)$ leave $W$ at $Out(\ww)$ at:
\begin{equation}
\begin{array}{c}
\Phi_{\ww }(r,\varphi )=(c_{3}-g_\ww\ln
r+\varphi,c_{4}r^{\delta_\ww})=(x,y)
\end{array},
\end{equation}
where $\delta_\ww$ is the \emph{saddle index} of $\ww$,
$$
\delta_\ww=\frac{C_{\ww }}{E_{\ww}}>0,
\quad
g_\ww=-\frac{\alpha_{\ww }}{E_{\ww }}<0,
\quad
c_{3}=g_\ww\ln \varepsilon
\quad
\text{ and }
\quad
c_{4}=\varepsilon ^{1-\delta_\ww }>0.
$$
The  minus sign in the equation $\dot{\theta}= - \alpha_{\ww }$ of  \eqref{localMapv} corresponds to the hypothesis (P\ref{H4}) and this implies that
$g_\ww<0$.

\subsection{Geometry near the cycle}
\label{Geom_vw}
\begin{figure}[hhh]
\begin{center}
\includegraphics[width=16cm]{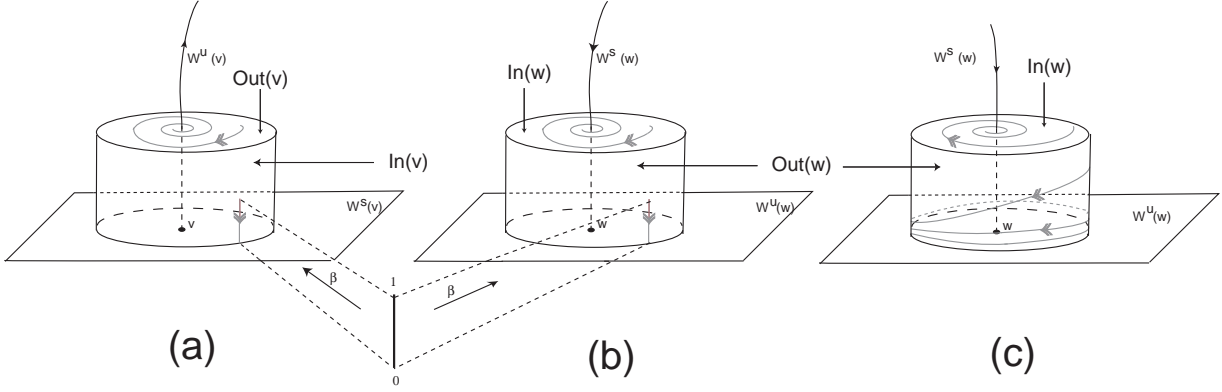}
\end{center}
\caption{\small Smooth structures referred in Lemma~\ref{Structures}.
The double arrows on the gray curves
(segment, spiral and helix) indicate correspondence of orientation and not the flow.
 (a) A segment $\beta $ in $In(\vv)$ is mapped by $\Phi _{\vv}$ into
 a spiral in $Out(\vv)$ around $W^u(\vv)$. (b) A segment $\beta $ in $Out(\ww)$ is mapped by $\Phi _{\ww}^{-1}$ into
 a spiral in $In(\ww)$ around $W^s(\ww)$. (c) If (P\ref{H4}) does not hold, a spiral in $In(\ww)$ around $W^s(\ww)$ is mapped
by $\Phi _{\ww}$ into a helix in $Out(\ww)$ accumulating on the circle  $Out(\ww)\cap W^{u}(\ww)$.
For the behaviour when (P\ref{H4}) holds see Figure~\ref{tangencies}.
}
\label{topological_structures}
\end{figure}

The notation and constructions of this section may now be used to describe the geometry associated to the local dynamics around the cycle. We start with some definitions that help the geometric description near each saddle-focus, illustrated in Figure \ref{topological_structures}.

\begin{definition}

\begin{enumerate}
\item \label{def_segmento} A \emph{ segment} $\beta $
 in  $In(\vv)$ or $Out(\ww)$ is a smooth regular parametrised curve of the type
$$\beta :(0,1]\rightarrow In(\vv) \qquad \text{or} \qquad \beta :(0,1]\rightarrow Out(\ww)$$ that meets $W^{s}_{loc}(\vv )$ or $W^{u}_{loc}(\ww )$ transversely at the point $\beta (0)$ only and such that, writing $\beta (s)=(x(s),y(s))$,
both $x$ and $y$ are monotonic and bounded functions of $s$ and $\frac{dx}{ds}$ is bounded.

\item \label{def_spiral}
A \emph{spiral} in $Out(\vv)$ or $In(\ww)$  around  a point $p$ is a curve  {without self-intersections}
$$\alpha :(0,1]\rightarrow Out(\vv) \qquad \text{or} \qquad \alpha :(0,1]\rightarrow In(\ww)$$
satisfying $\dpt \lim_{s\to 0^+}\alpha (s)=p$ and such that, if
$\alpha (s)=(\alpha _{1}(s),\alpha _{2}(s))$ are its expressions  in
polar coordinates $(\rho ,\theta )$ around $p$, then
 {$\alpha_2$ is monotonic with $\lim_{s\to 0^+}|\alpha _{2}(s)|=+\infty$.}

\item  Consider a cylinder $C$ parametrised by a covering $(\theta,h )\in  \RR\times[a,b]$,
with $a<b\in\RR$ where $\theta $ is periodic.
A \emph{helix} in the cylinder $C$
\emph{accumulating on the circle}
$h=h_{0}$ is a curve
$\gamma :(0,1]\rightarrow C$   {without self-intersections}
such that its coordinates $(\theta (s),h(s))$
satisfy
$$
\lim_{s\to 0^+}h(s)=h_{0}
\qquad\qquad
\lim_{s\to 0^+}|\theta (s)|=+\infty
$$
 {with $h$ monotonic.}
%
\end{enumerate}
\end{definition}

 We are interested in spirals for which the  point $p$ in  the definition is  the intersection of the
 two-dimensional local stable/unstable manifold of $\vv$ or $\ww$ with the  cross section.
The next lemma from Aguiar \emph{et al} \cite[\S 6]{ALR} summarises some basic
 results on the geometry near the saddle-foci.
In its original form the authors assume implicitly that  property (P\ref{H4}) does not hold,
 {but the same proof holds, with  minor adaptations.}

\begin{lemma}[Aguiar \emph{et al} \cite{ALR}]
\label{Structures}
A segment $\beta $:
\begin{enumerate}
\item
in $In(\vv)$ is mapped by $\Phi _{\vv}$ into  a spiral in $Out(\vv)$ around $W^u(\vv)$;
\item
in $Out(\ww)$ is mapped by $\Phi _{\ww}^{-1}$ into a spiral in $In(\ww)$ around $W^s(\ww)$;
\item \label{spiral3}
in $In(\vv)$ is mapped by $\Phi _{\vv}$ into  a spiral in $Out(\vv)$ around $W^u(\vv)$,
that is mapped by $ \Psi_{\vv \rightarrow \ww}$  into another spiral around $W^s(\ww)\cap In(\vv)$.
 \end{enumerate}
\end{lemma}

 If (P\ref{H4}) does not hold, the spiral in \ref{spiral3}. is mapped by $\Phi _{\ww}$ into a helix in $Out(\ww)$ accumulating on the circle  $Out(\ww) \cap W^{u}(\ww)$. 
 A key argument in obtaining heteroclinic tangencies is that  this fails when   (P\ref{H4}) holds.

\begin{figure}[hhh]
\begin{center}
\includegraphics[height=3cm]{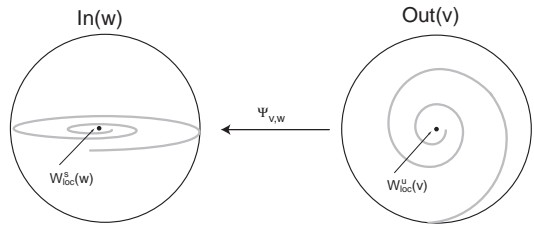}
\end{center}
\caption{\small The transition map from $V$ to $W$ may be approximated by
$\Psi_{\vv, \ww}(x,y)=\left(ax, y/a\right) $, where $(x,y)$ are the rectangular coordinates at $Out(\vv)$ and $In(\ww)$.
A circle with centre at $W^s_{loc}(\vv)$ is mapped into an  ellipse and a spiral is deformed in a similar way.}
\label{spiral_shrink}
\end{figure}

The transition map $\Psi_{\vv,\ww}$ has a simple geometry,  {shown} in Figure \ref{spiral_shrink}.

\begin{lemma}
A circle of radius $r<\varepsilon$ in $Out(\vv)$ centered at the origin is mapped by $\Psi_{\vv,\ww}$ into an ellipse centered at the origin of $In(\ww)$ with major axis of length $ar\ge r$ and  minor axis of length $\frac{r}{a}\le r$.
\end{lemma}
Note that the map $\Psi_{\vv,\ww}$ is given in rectangular coordinates. To compose this map with $\Phi _{\ww }$, it is required to change the coordinates. We address this issue in Section \ref{SectionT-pointHetero}.

\section{Heteroclinic Tangencies and Horseshoes}

\label{SectionT-pointHetero}

In this section we give a more precise formulation of Theorems \ref{teoremaReverteRotacao} and \ref{TeoremaTransversal} and prove them.
 In order to simplify computations, we assume from now on, that $W^u(\ww) \cap In(\vv)$ and $W^s(\vv) \cap Out(\ww)$ are vertical segments across $In(\vv)$ and $Out(\ww)$, respectively.
Let $\beta(s)=(0,s) \subset In(\vv)$, $s \in [0,\varepsilon]$, be a parametrisation of $W^u(\ww) \cap In(\vv)$, where $(0,0)$ are the local coordinates of the point $[\vv \rightarrow \ww] \cap In(\vv)$.
Then $$\Phi_{\ww} \circ \Psi_{\vv ,\ww  } \circ \Phi_{\vv}(\beta(s))$$ defines an oriented curve in $Out(\ww)$. 
Our first step is to show that, under hypotheses (P\ref{H1})--(P\ref{H4}), 
 this curve changes the direction of its turning around the cylinder $Out(\ww)$ at infinitely many points where it has a vertical tangent,  {when the parameters that determine the linear part of the vector field at $\vv$ and $\ww$ lie in a set of full Lebesgue measure.}

\subsection{Preliminary analysis of maps}
\label{preliminary_analysis}
 {We start by obtaining}
 the expression for $\Phi_{\ww} \circ \Psi_{\vv ,\ww } \circ \Phi_{\vv }(\beta(s))$. The following result is a generalisation of Ovsyannikov and Shilnikov \cite{OS}, treating the dependence on $a$ and containing simpler expressions.

\begin{lemma}
\label{x,y}
Let $\beta=(0, s)$, $s \in [0,\varepsilon]$, be a segment in $In(\vv)$ parametrized by $s$ where $$\{(0, 0)\} = [\ww \rightarrow \vv] \cap In(\vv)$$ and let $(x_\ww(s),y_\ww(s))= \Phi_{\ww} \circ \Psi_{\vv ,\ww } \circ \Phi_{\vv } (\beta(s))$. Then:

\begin{equation}\label{xweyw}
\left\{
\begin{array}{l}
x_\ww(s)=
-g_\ww \delta_\vv\ln s - \frac{g_\ww}{2} \ln C(\varphi)+ \Phi(\varphi) + c_3 - g_\ww \ln c_1\\
 \\
y_\ww(s) =
c_4 c_1^{\delta_\vv \delta_\ww} s^{\delta_\vv \delta_\ww} [C(\varphi)]^{\frac{\delta_\ww}{2}}\\
\end{array}
\right.
\end{equation}
where
\begin{equation}\label{CePhi}
\left.
\begin{array}{l}
\varphi(s)=-g_\vv \ln s +c_2,\\
\\
C(\varphi)= a^2 \cos^2(\varphi) + \frac{1}{a^2}  \sin^2(\varphi) \quad \text{and}\\
\\
\Phi(\varphi)= \arg \left(a \cos(\varphi),\frac{1}{a} \sin(\varphi) \right),\\
\end{array}
\right.
\end{equation}
with the argument $\arg$  taken in the  interval $\left[\frac{k\pi}{2}, \frac{(k+1)\pi}{2}\right]$, $k \in \mathbf{Z}$ that contains $\varphi$.
\end{lemma}

\textbf{Proof:}
The image of $\beta$ under the local map near $\mathbf{\vv}$, $\Phi_{\vv}$, is given by:
\begin{equation}
\label{closesym1}
\Phi_{\vv}(0,s) = (c_1 s^{\delta_\vv}, -g_\vv \ln s + c_2)=(r, \varphi)
\end{equation}
that in rectangular coordinates has the form $(X(s), Y(s))=(r(s)\cos \varphi(s), r(s)\sin \varphi(s))$, hence:

\begin{equation}
\label{local map vw}
\Psi_{\vv ,\ww  }\circ \Phi_{\vv}(0,s)= \left(a\  r(s) \cos \varphi(s), \frac{r(s)}{a}\sin \varphi(s)\right)
\end{equation}

Therefore the radial component of $\Psi_{\vv ,\ww } \circ \Phi_{\vv } (\beta(s))$ may be written as $R=r(s) \sqrt{C(\varphi)}$ and the angular component is
$$
\Phi= \arg \left( c_1 s^{\delta_\vv}\left(a\cos(\varphi),\frac{1}{a} \sin(\varphi)\right) \right)=
 \arg \left( a\cos(\varphi),\frac{1}{a} \sin(\varphi)\right)
$$
since $c_1 s^{\delta_\vv}>0$.
Using the expression $\Phi_{\ww }$ of the local map near $\mathbf{w}$ it follows that $\Phi_{\ww} \circ \Psi_{\vv ,\ww } \circ \Phi_{\vv } (\beta(s))$ has the form given in the statement of the lemma.
\Qed

We can now give a description of the global dynamics near the whole network. Hereafter, denote by
$\delta$ the product $\delta_\vv \delta_\ww$,  and by $\eta$ the map $\Phi_{\ww }\circ \Psi_{\vv,\ww }\circ \Phi_{\vv }: In(\vv) \rightarrow Out(\ww)$,  {and} let $\dpt\gamma=\frac{\alpha_\ww}{\alpha_\vv}\frac{C_\vv}{E_\ww}$.

\begin{lemma}
\label{Lemma_aux2}
Let $\beta$ be a segment on $In(\vv)$ parametrized in rectangular coordinates by $(x_\vv(s), y_\vv(s))$ and let $(x_\ww(s), y_\ww(s))$ be the coordinates of $\eta(\beta(s))\in Out(\ww)$. Then:
\begin{enumerate}
\item\label{Case1} if $a = 1$, then the coordinates $x_\ww(s)$ and $y_\ww(s)$ are both monotonic functions of $s$;
\item\label{Case2}  if $a \neq 1$, then the coordinate $y_\ww(s)$ is not a monotonic function of $s$;
\item\label{Case3} for all $a\geq 1$, $\lim_{s \rightarrow 0^+} y_\ww(s)=0$;
\item\label{newCase}
 { for all $a\geq 1$, }
if $\gamma> 1$ then $\lim_{s \rightarrow 0^+} x_\ww(s)=-\infty$ and
if $0<\gamma<1$ then $\lim_{s \rightarrow 0^+} x_\ww(s)=+\infty$.
\end{enumerate}
\end{lemma}

\textbf{Proof:}

{\sl \ref{Case1}.}\
If $a=1$, then $C(\varphi)=a^2 \cos^2 (\varphi) + \frac{1}{a^2} \sin^2(\varphi) \equiv 1$.
Hence, $y_\ww(s)=k_1 s^{\delta}$, where $k_1 \in \RR^+$ and the result follows directly.
Since $a=1$, we may write
$x_\ww(s)= k_2-(g_\ww \delta_\vv + g_\vv) \ln s=k_2-\frac{\alpha_\vv}{E_\vv}\left(1-\gamma\right)\ln s$,
where $k_2 \in \RR^+$, then $x_\ww(s)$ is monotonic.
If $\gamma\ne 1$ then $x_\ww(s)$ is strictly monotonic.

{\sl \ref{Case2}.}\
If $a \neq 1$, then $y_\ww(s)$ is not monotonic because, in $In(\ww)$, the euclidean distance  between $\Psi_{\vv,\ww }\circ \Phi_{\vv }(\beta(s))$ and $W^s(\ww)$ is not a decreasing function of the parameter $s$. Thus, the map $y_\ww(s)$ that represents the height is not monotonic. Recall that $\Psi_{\vv,\ww }$ consists of an expansion in the horizontal direction and a contraction in the vertical direction.

{\sl \ref{Case3}.}\
Since $(C(\varphi))^\frac{\delta_\ww}{2}=\left(a^2 \cos^2 (\varphi) + \frac{1}{a^2} \sin^2(\varphi)\right)^\frac{\delta_\ww}{2} $ is bounded, then $\displaystyle\lim_{s \rightarrow +\infty} y_\ww(s)= \lim_{s \rightarrow +\infty} k_1 s^{\delta_\vv}=0.
$

{\sl \ref{newCase}.}\
Let $\hat{x}(s)=-g_\ww \delta_\vv\ln s + \Phi(\varphi(s)) $.
Note that $\hat{x}(s)-x_\ww(s)=- \frac{g_\ww}{2} \ln C(\varphi)+ c_3 - g_\ww \ln c_1$ is limited, hence it is sufficient to compute $\lim_{s \rightarrow 0^+}\hat{x}(s)$.

Since $\lim_{s \rightarrow 0^+}\varphi(s)=+\infty$, then as ${s \rightarrow 0^+}$ one gets
$\dpt\frac{2}{\pi}\varphi(s)\in\left[ k, k+1\right]$
with $k\in{\rm\bf N}$,  $k\rightarrow +\infty$ and thus
$\dpt\frac{2}{\pi}\Phi(s)$ lies in the same interval.
Using $\dpt-g_\ww \delta_\vv\ln s=\frac{g_\ww\delta_\vv}{g_\vv}\varphi(s)=-\gamma\varphi(s)$
it follows that
$\dpt\frac{2}{\pi}\hat{x}(s)\in\left[(1-\gamma)k-\gamma,(1-\gamma)k+1 \right].$
Hence, if $0<\gamma<1$ then $\dpt\lim_{s \rightarrow 0^+}\hat{x}(s)=+\infty$, and if $\gamma>1$ then
$\dpt\lim_{s \rightarrow 0^+}\hat{x}(s)=-\infty$.
\Qed

The dynamics in the reversible case, when $\gamma=1$, has been described in \cite{KLW, LTW}.

\subsection{The reversal property}
\label{The reversal property}
\begin{figure}
\begin{center}
\includegraphics[height=8cm]{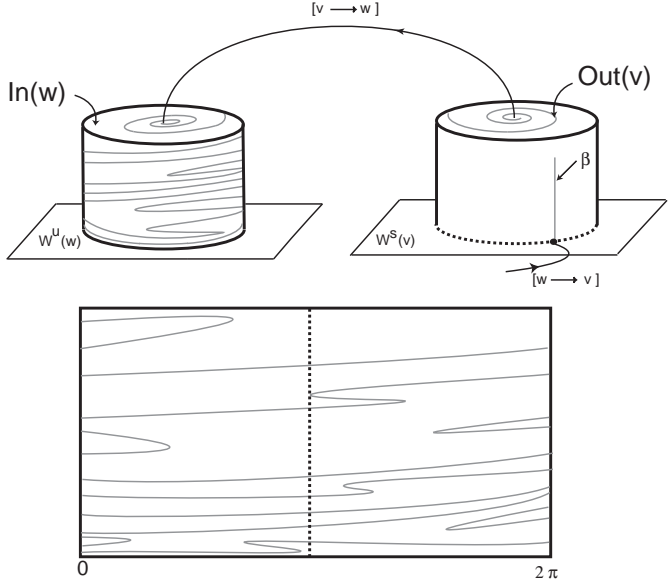}
\end{center}
\caption{\small Top:
the vertical segment $\beta$ on $In(\vv)$ is mapped by $ \Psi_{\vv,\ww }\circ \Phi_{\vv }$ into a distorted spiral in $In(\ww)$ and by $\eta$ into a curve in $Out(\ww)$ that accumulates on
$W^u(\ww)$.
Bottom:
on the cross section $Out(\ww)$,  the curve $\eta\circ \beta$ reverses the orientation of its angular coordinate infinitely many times as it accumulates on $W^u(\ww)$ and crosses $W^s(\vv)$ (dotted line) infinitely many times. The angular coordinates of the points of reversal are dense in the circle,  small changes in $W^s(\vv)\cap Out(\ww)$ create tangencies that
coexist with transverse crossings of the invariant manifolds.}
\label{tangencies}
\end{figure}

The main goal now is to prove that under hypothesis (P\ref{H1})--(P\ref{H4}), the coordinate map
$x_\ww$ is not a monotonic  function of $s$, since the curve $\eta \circ \beta$ reverses the direction of its turning around $Out(\ww)$ infinitely many times. This is the notion illustrated in Figure \ref{tangencies} and formalized in the following definition:

\begin{definition}
We say that the vector field $f$ has the \emph{dense reversals property} if for the vertical segment $\beta(s)=(0, s)\in In(\vv)$, $s \in [0, \varepsilon]$, the projection into $W^u_{loc}(\ww)$ of the points where $\eta\circ \beta$ has a vertical tangent is dense in $W_{loc}^u(\ww)\cap Out(\ww)$.
\end{definition}

The dense reversals property is the key step in the proof of Theorem \ref{teoremaReverteRotacao}.
In order to prove it we need some additional assumptions on  the parameters
$P=\left(\alpha_\vv, C_\vv, E_\vv, \alpha_\ww, C_\ww, E_\ww\right)$,
that determine the linear part of the vector field $f$ at the equilibria.
Note that Case \ref{Case1} of Lemma \ref{Lemma_aux2} rules out reversals when $a=1$.
For any  $a>1$, let $\mathcal{B}$ be the subset of parameters  given by:
\begin{equation}\label{final_formulae}
\mathcal{B}=\left\{ P:\
\left(a^2-\frac{1}{a^2}\right)\frac{2\alpha_\vv}{C_\vv - \sqrt{\alpha_\vv^2+ 4C_\vv^2}}
\leq \frac{E_\ww}{\alpha_\ww}- \frac{a^2 C_\vv}{\alpha_\vv} \leq
\left(a^2-\frac{1}{a^2}\right)\frac{2\alpha_\vv}{C_\vv + \sqrt{\alpha_\vv^2+4 C_\vv^2}}
\right\}
\end{equation}
and having non empty interior $int(\mathcal{B})$.
Let $\mathcal{D}$ be the dense subset  of $\mathcal{B}$ given by:
\begin{equation}\label{denseSet}
\mathcal{D}=\left\{ P\in int(\mathcal{B}):\
\gamma=\frac{\alpha_\ww}{\alpha_\vv}\frac{C_\vv}{E_\ww} \notin \textbf{Q}
\right\} .
\end{equation}

The condition (\ref{final_formulae}) in the definition of  $\mathcal{B}$ is satisfied by an open set of parameters
$\alpha_\vv$, $C_\vv$, $\alpha_\ww$, $E_\ww$ and does not involve the quantities $E_\vv$ and $C_\ww$.
The  condition \eqref{denseSet} defining $\mathcal{D}$ also does not
 involve $E_\vv$  and $C_\ww$. 
The next result shows that  the dense reversals property is generic in the class of  vector fields that satisfy (P\ref{H1})--(P\ref{H4}) with parameters in   $\mathcal{B}$.
This will be used to complete the proof Theorem~\ref{teoremaReverteRotacao}.

\begin{theorem}
\label{prop12}
Let $f$ be a vector field on $\EU^3$ satisfying (P\ref{H1})--(P\ref{H4}). If the parameters for the linear part of $f$ near $\vv$ and $\ww$ lie in $\mathcal{D}$  with $a>1$ then
$f$ has the \emph{dense reversals property}.
\end{theorem}

\begin{figure}[h]
\begin{center}
\includegraphics[width=10cm]{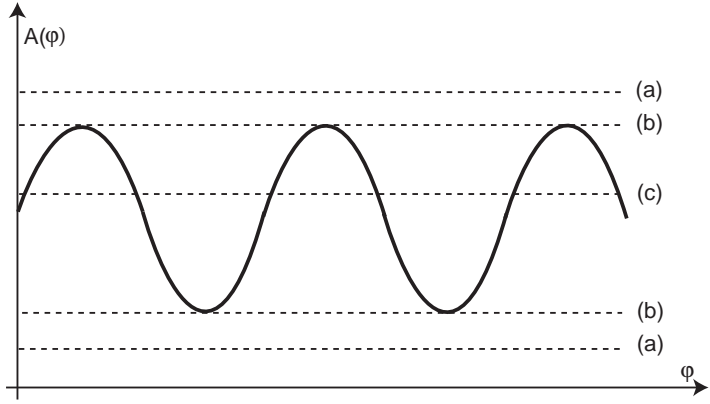}
\end{center}
\caption{\small Positions of the graph of $A(\varphi)$ of \eqref{AdePhi} (solid curve) with respect to different values of
$\alpha_\vv E_\ww/\alpha_\ww$ (dashed lines).
Parameters outside the set $\mathcal{B}$ of \eqref{final_formulae} correspond to cases (a), where the graph of $A(\varphi)$ never crosses the dashed line, and hence  $x_\ww(s)$ is a monotonic function of $s$.
Cases (b) and (c) correspond to parameters  in $\mathcal{B}$.
Points where  $x_\ww(s)$ has a vertical tangent are created in case (b) as  the graph of
$A(\varphi)$ is tangent to the dashed line, for parameters in the boundary of  $\mathcal{B}$.
The graph of
$A(\varphi)$ crosses the dashed line transversely in case (c) for parameters in the interior of
  $\mathcal{B}$, giving rise to two $\pi$-periodic sequences of points  of vertical tangency.
}
\label{plots}
\end{figure}

\textbf{Proof: }
We need to compute the coordinate $x_\ww(s)$ at the points where $\eta(\beta(s))$ has a vertical tangent, where $\beta(s)=(0,s) \subset In(\vv)$, $s \in (0,\varepsilon]$ is a parametrisation of a vertical segment. Differentiating the expression \eqref{xweyw} of Lemma \ref{x,y}, we get:
\begin{equation}\label{eqdxw}
\frac{dx_\ww}{ds}=-\frac{1}{s}\left[g_\ww\delta_\vv + \frac{1}{C(\varphi)}
 \left(2g_\ww g_\vv \left(a^2-\frac{1}{a^2}\right)\sin\varphi \cos\varphi + g_\vv\right)\right]
\end{equation}
and hence ${dx_\ww}/{ds}=0$ has solutions if and only if
$A(\varphi)={\alpha_\vv E_\ww}/{\alpha_\ww}$ where
\begin{equation}\label{AdePhi}
A(\varphi)=
C_\vv a^2 \cos^2\varphi + \frac{C_\vv}{a^2}\sin^2\varphi +
2 \alpha_\vv \left(a^2-\frac{1}{a^2}\right)\sin\varphi\cos\varphi.
\end{equation}
Therefore ${dx_\ww}/{ds}=0$ has solutions (see Figure~\ref{plots}) if and only if
\begin{equation}
\label{max}
\min A(\varphi) \leq \frac{\alpha_\vv E_\ww}{\alpha_\ww} \leq \max A(\varphi).
\end{equation}

In order to determine the maxima and the minima of $A(\varphi)$, we compute:
$$
\frac{dA}{d\varphi}=2 \left(a^2-\frac{1}{a^2}\right) \left[\alpha_\vv \cos^2\varphi -
C_\vv \cos\varphi\sin\varphi-\alpha_\vv\sin^2\varphi \right]
=2\left(a^2-\frac{1}{a^2}\right)\mathcal{Q}(\cos\varphi,\sin\varphi),
$$
where $\mathcal{Q}(x,y)=\alpha_\vv x^2 - C_\vv xy - \alpha_\vv y^2$.
The quadratic form $\mathcal{Q}(x,y)$ equals zero on the lines that join the origin to the points
$$
(x_\star^\pm, y_\star)=(C_\vv \pm \sqrt{C_\vv^2+4\alpha_\vv^2},2 \alpha_\vv).
$$
Hence
$$
\frac{dA}{d\varphi}=0\quad\Leftrightarrow\quad
(\cos\varphi_\star^\pm,\sin\varphi_\star^\pm)= \frac{1}{N}\left(x_\star^\pm, y_\star\right) \qquad \text{with} \qquad N=\|(x_\star^\pm, y_\star)\|.
$$
In order to impose condition (\ref{max}), we write $A(\varphi)=\mathcal{R}(\cos\varphi,\sin\varphi)$ where
$$
\mathcal{R}(x,y)= C_\vv a^2 x^2 + 2\alpha_\vv \left(a^2 - \frac{1}{a^2}\right)xy + \frac{C_\vv}{a^2}y^2.
$$
Since $\mathcal{R}(x,y)>0$ if $x>0$ and $y>0$, then $(x_\star^+, y_\star)$ corresponds to maxima of $A(\varphi)$ and $(x_\star^-, y_\star)$ to minima.

Let $S=\sqrt{C_\vv+4\alpha_\vv^2 }$, then
$$
N^2=2C_\vv^2+8\alpha_\vv^2\pm 2C_\vv S=2 S\left(S\pm C_\vv\right)$$
and
$$
\mathcal{R}(x_\star^\pm, y_\star)=C_\vv a^2 N^2\pm 2\alpha_\vv^2 \left(a^2 - \frac{1}{a^2}\right) S
$$
and hence
$$
A(\varphi^\pm_\star)= \frac{1}{N^2} \mathcal{R}(x_\star^\pm, y_\star)=
a^2 C_\vv \pm 2 \alpha_\vv^2 \left(a^2 - \frac{1}{a^2}\right)\frac{S}{2 S\left(S\pm C_\vv\right)}=
a^2 C_\vv + 2 \alpha_\vv^2 \left(a^2 - \frac{1}{a^2}\right)\frac{1}{C_\vv\pm S},
$$
showing that conditions (\ref{final_formulae}) and (\ref{max}) are equivalent.
These conditions imply that $A(\varphi)=\alpha_\vv E_\ww/\alpha_\ww$ at infinitely many values
$\varphi=\varphi_0+n\pi,$
where $\varphi_0 \in [0, \pi]$ and $n \in \textbf{Z}$. Since $\varphi=-g_\vv \ln s + c_2$, then
${dx_\ww}/{ds}=0$ has solutions
$$s_n=s_0e^{- \frac{n\pi}{g_\vv}}, \quad n=0, 1,2, \ldots \qquad \text{where} \qquad s_0= e^{- \frac{\varphi_0}{g_\vv}}e^{\frac{c_2}{g_\vv}}.
$$
In Lemma~\ref{lemaxwsn} below, we show that  $x_\ww\left(s_n\right)=x_\ww(s_0)+n\pi \left(1-\gamma\right)$.
Hence, if the genericity condition $\gamma\notin\textbf{Q}$ in \eqref{denseSet} holds, then the points $x_\ww(s_n)$ are dense in the circle defined by $W^u_{loc}(\ww) \cap Out(\ww)$.
\Qed

By the implicit function theorem, the arguments that we have used for vertical tangencies are still valid if $\beta(s)$ is  any line with slope close to the vertical.

\begin{lemma}\label{lemaxwsn}
For any $s_0\in\RR$ and $n=0, 1,2, \ldots$, we have
\begin{equation}\label{xwsn}
x_\ww\left(s_0 e^{\frac{-n\pi}{g_\vv}}\right)=x_\ww(s_0)+n\pi \left(1-\gamma\right)
\qquad\mbox{for}\qquad
\gamma=\frac{\alpha_\ww}{\alpha_\vv}\frac{C_\vv}{E_\ww} .
\end{equation}
\end{lemma}

\textbf{Proof:}
Using the expressions \eqref{CePhi} we get that $\varphi(s)=\varphi(s_0)+n\pi$ and
 $C(\varphi+\pi)=C(\varphi)$.
Also, if $\varphi(s) \in \left[\frac{k\pi}{2},\frac{(k+1)\pi}{2}\right]$, then $\Phi(\varphi)$ lies in the same interval and
$\Phi(\varphi+\pi)=\Phi(\varphi)+\pi$.
The result follows, using the expression \eqref{xweyw}  for $x_\ww(s)$ in Lemma \ref{x,y}.
\Qed

The next theorem is a precise formulation of Theorem \ref{teoremaReverteRotacao} and justifies the title of the article.

\begin{theorem}
\label{teoremaFinal}
Let $\mathcal{A}_P$ be the set of  vector fields on $\EU^3$ satisfying (P\ref{H1})--(P\ref{H4}), for a given parameter  value $P\in \mathcal{B}$.
Then there is an open subset $\mathcal{C}$ of  $\mathcal{A}_P$ such that the
 set of vector fields in $\mathcal{C}$ whose flow has
tangencies between the two dimensional invariant manifolds $W^u(\ww)$ and $W^s(\vv)$ is dense in $\mathcal{A}_P$ in the $C^k$ topology, for every $k\ge 2 \in \mathbf{N}$.
 \end{theorem}

 \textbf{Proof:}
 {Suppose} that $W^u_{loc}(\ww) \cap In(\vv)$ and $W^s_{loc}(\vv) \cap Out(\ww)$ are vertical segments across $In(\vv)$ and $Out(\ww)$, respectively.
Let $\beta(s)=(0,s) \subset In(\vv)$, $s \in (0,\varepsilon]$, be a parametrisation of $W^u_{loc}(\ww) \cap In(\vv)$, where $(0,0)$ is the point $[\vv \rightarrow \ww] \cap In(\vv)$ and
let $\gamma(s)=(x_0,s)$, $s\in(0,\varepsilon]$ be a parametrisation of
$W^s_{loc}(\vv)\cap Out(\ww)$, where  {$(x_0,0)$} is the point $[\vv\to\ww]\cap Out(\ww)$.

If $x_0$ is  the projection of a point where
$\eta\circ \beta$ has a vertical tangent, then  $W^s(\vv)$ is tangent to $W^u(\ww)$ at the corresponding point. Otherwise, by Theorem~\ref{prop12},  {there is a point $(x_1,y_1)$ in $\eta\circ\beta$ with $|x_1-x_0|$ as small as we want and $y_1>0$.  A small modification of the map  $\Psi_{\ww,\vv}$ around this point will not affect the curve $\eta\circ\beta$ but will move $W^s_{loc}(\vv)\cap Out(\ww)$,}
creating the tangency.

 {The open set $\mathcal{C}$ in the statement of Theorem~\ref{teoremaFinal} consists of vector fields for which $W^s_{loc}(\vv)\cap Out(\ww)$ is close to a vertical segment.}
Theorem~\ref{prop12} still holds if  $W^u_{loc}(\ww) \cap In(\vv)$ is a line with slope close to the vertical,
 {as remarked at the end of its proof.}
Hence  {in this case} we still have the dense reversals property   {and therefore we obtain Theorem~\ref{teoremaFinal} } as we proceed to explain.
Suppose  $W^s_{loc}(\vv)\cap Out(\ww)$ is a curve close to a vertical segment and that it is parametrised by
$\xi(s)=(x(s),y(s))$, with $\xi(0)=(x_0,0)$.
Then,
 {if necessary, we may change $\Psi_{\ww,\vv}$ around a point $(x_1,y_1)$ as  above to have $\xi(s)$ meeting }
$\eta \circ \beta$ at a point where the last curve  has a vertical tangent.
Since $\xi$ is close to a vertical line, then its slope near the intersection is close to the vertical,
so a  {second} small change of the  {transition map}
near this point will create a tangency. 
\Qed

 {The proof of Theorem~\ref{teoremaFinal}}
allows us to connect the unstable manifold of $\ww$ with the stable manifold of $\vv$ without  {recourse to} the Pasting and Connecting Lemmas  {(see \cite{Hayashi})}. 
 {Using}
these results in $In(\vv)$ for the first return map,  {would} not guarantee
 {that the perturbed diffeomorphism would be  the first return map of a vector field with a Bykov cycle.}

\subsection{Topological and Hyperbolic Horseshoes}
\label{Topological and Hyperbolic Horseshoes}
In this section we give the geometrical construction for the proof of
Theorem~\ref{TeoremaTransversal}.  This is the standard construction for establishing symbolic dynamics, except for the obstacle of  infinitely many reversions, that we will overcome by showing that the non-transverse intersections can be avoided,
since generically the line $\eta\circ \beta(s)$ intersects transverselly $W^s_{loc}(\vv)\cap Out(\ww)$ infinitely many times. This phenomenon coexists with the denseness of the tangencies in $\mathcal{A}_P$.

First we need to recall some terminology about horizontal and vertical strips used in Guckenheimer and Holmes \cite{Guckenheimer and Holmes} adapted to our purposes.
Given a rectangle $[w_1,w_2]\times [z_1, z_2]$ in either $In(\vv)$ or  $Out(\ww)$,
 a \emph{horizontal strip} across the rectangle is the set
$$
\{(x,y): x \in [w_1,w_2], y\in[u_1(x),u_2(x)]\},
$$
where $ u_1,u_2: [w_1,w_2] \rightarrow [z_1,z_2]$  are Lipschitz functions such that $u_1(x)<u_2(x)$. The \emph{horizontal boundaries} of the strip are the graphs of the $u_i$ and
the \emph{vertical boundaries} are the lines  $\{w_i\}\times  [u_1(w_i),u_2(w_i)]$.
A \emph{vertical strip} across the rectangle has a similar definition with the roles of $x$ and $y$ reversed.

Recall that $P=(\alpha_\vv, C_\vv, E_\vv,\alpha_\ww, C_\ww, E_\ww)$ is the set of parameters that determine the linear part of $f$ at the hyperbolic saddle-foci
and that  $\beta(s)=(0, s)\in In(\vv)$, $s \in (0, \varepsilon]$.
The next result holds for almost all vector fields on $\EU^3$ satisfying (P\ref{H1})--(P\ref{H4}).
One exception is the case $\gamma=1$, that occurs in reversible vector fields, as remarked before.
The other  exception occurs when $\gamma\in\mathbf{Q}$ and $P\in\mathcal{B}$ and when moreover there are points in $Out(\ww)$ with first coordinate equal to zero, where $\eta\circ\beta$ has a vertical tangent --- we say in this case that the vector field has a \emph{periodic tangency}.

 \begin{proposition}\label{propStrips}
If $\displaystyle\gamma=\frac{\alpha_\ww}{\alpha_\vv}\frac{C_\vv}{E_\ww} \ne 1$, then
for all vector fields on $\EU^3$ satisfying (P\ref{H1})--(P\ref{H4}) and not having a periodic tangency, the following holds: for any sufficiently small $\tau>0$, there is a sequence of disjoint horizontal strips across the rectangle $[0,\tau]\times[0,\tau]\subset In(\vv)$, accumulating on $W^s(\vv)$, whose image by the first return map $\Psi_{\ww,\vv}\circ \eta$ is a vertical strip across $[0,\tau]\times[0,\tau]$.
 \end{proposition}

\textbf{Proof:}
Take  $(0,0)$ as the local coordinates of the points $[\vv \rightarrow \ww] \cap In(\vv)$ and
$[\vv \rightarrow \ww] \cap Out(\ww)$, as before.
Start with the rectangle $\left[0,\tau_0\right]\times\left[0,\tau_0\right]\subset In(\vv)$ with
$0<\tau_0<\min\{\pi,\varepsilon\}$.
For each $t \in \left[0,\tau_0\right]$, define the family of vertical segments
$\beta_t(s) = (t, s)\in In(\vv)$, and let $\beta(s)=\beta_0(s)$ with  $x_\ww(s)$ the first coordinate of $\eta(\beta(s))$, as before.
Then  $\eta(\beta_t(s))$, for different $t$, are disjoint curves in $Out(\ww)$.
We will assume from now on that $\gamma>1$, the proof in the case $0<\gamma<1$ is obtained by replacing increasing functions by decreasing functions, and $-\infty$ by $+\infty$ in what follows.

The proof consists in finding $\tau$ with $0<\tau\le\tau_0$ and a decreasing sequence $s_n>0$, with
$\lim_{n\to\infty} s_n=0$,
such that in each interval $\left(s_{2n+1},s_{2n}\right)$ the function $x_\ww(s)$ is monotonically increasing  and crosses $\left[-\tau,0\right]\pmod{2\pi}$. Hence on this interval the curve
$\eta(\beta(s))$ is transverse to each vertical line in $Out(\ww)$ and for each $n$ there are
$a_n<b_n\in\left(s_{2n+1},s_{2n}\right)$ such that $x_\ww(a_n)=-\tau\pmod{2\pi}$ and $x_\ww(b_n)=0\pmod{2\pi}$.

The curves $\eta(\beta_t(s))$ have the same properties, i.e., by taking a smaller $\tau>0$ if necessary, since the curves depends smoothly on $t$, the following holds:
for each $t\in\left[0,\tau\right]$ there are two sequences $0<a_n(t)<b_n(t)<a_{n-1}(t)$, with
$\lim_{n\to\infty} a_n(t)=0$, such that in each interval $\left[a_n(t),b_n(t)\right]$ the first coordinate of
$\eta(\beta_t(s))$ is a monotonically increasing function of $s$, taking the values $-\tau\pmod{2\pi}$ at $s=a_n(t)$ and $0\pmod{2\pi}$ at $b_n(t)$.
This means that the strip across $\left[0,\tau\right]\times\left[0,\tau\right]\subset In(\vv)$ with horizontal boundaries $a_n(t)$ and $b_n(t)$ is mapped by $\eta$ into a horizontal strip across
$\left[-\tau,0\right]\times\left[0,\tau\right]\subset Out(\ww)$, that  in turn is mapped by $\Psi_{\ww,\vv}$ into a vertical strip  across $\left[0,\tau\right]\times\left[0,\tau\right]\subset In(\vv)$, as required.

\begin{figure}
\begin{center}
\includegraphics[width=15cm]{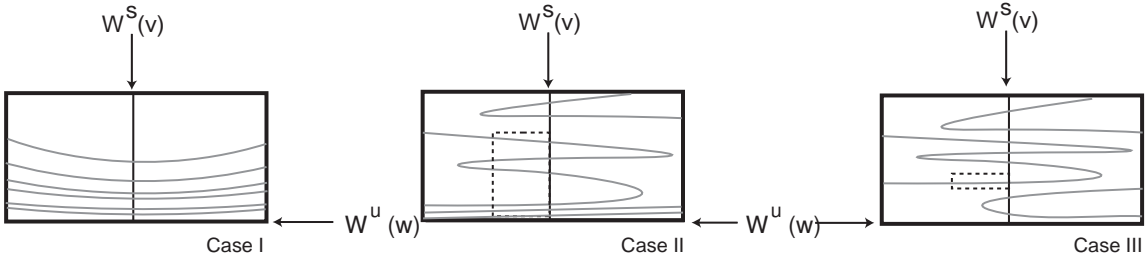}
\end{center}
\caption{\small The curve  $\eta(\beta(s))\subset Out(\ww)$ in the three first cases of the proof of Proposition~\ref{propStrips}.
 \textbf{Case I:} If $P \notin \mathcal{B}$, then $x_\ww(s)$ is monotonic.
 \textbf{Case II:} If $P  \in int( \mathcal{B} )\backslash \mathcal{D}$ with  $\gamma\in\textbf{Q}$, then reversion points are periodic.
  \textbf{Case III:} If $P \in \mathcal{D}$, then reversion points are dense.
  In cases \textbf{II} and \textbf{III} we restrict the strip to the dotted rectangle to avoid reversion points.}
\label{horseshoes}
\end{figure}

In order to obtain the sequence of monotonicity intervals we  distinguish four  cases (see Figure~\ref{horseshoes}):
\medbreak

\textbf{Case I:} If $P \notin \mathcal{B}$, then the curve $\eta(\beta(s))$ does not reverse the direction of its turning around $Out(\ww)$ and  $x_\ww$ is a monotonically  increasing function of $s$,  since $\gamma>1$.  Using Lemma~\ref{Lemma_aux2} it follows that $x_\ww(s)$ goes across $\left[-\tau,0\right]\pmod{2\pi}$ infinitely many times, as required.
\medbreak

\textbf{Case II:} If $P  \in int( \mathcal{B} )\backslash \mathcal{D}$, then the projection into $W^u_{loc}(\ww)\cap Out(\ww)$ of the points where the curve $\eta(\beta(s))\in Out(\ww)$ reverses orientation is finite.
Since we are assuming that the vector field does not have a periodic tangency, the curve
$\eta(\beta(s))$ is never tangent to the segment $(0,s)\subset Out(\ww)$.
Then there is a $\tau>0$ such that $dx_\ww/ds$ is never zero when
$x_\ww(s)\in\left[-\tau,0\right]\pmod{2\pi}$ and therefore $x_\ww(s)$ is a monotonic function of $s$ when $x_\ww(s)$ lies in that interval.
Since $\lim_{s\to 0^+} x_\ww(s)=-\infty$, it follows that $x_\ww(s)$ crosses $\left[-\tau,0\right]\pmod{2\pi}$ infinitely many times, as a monotonically increasing function of $s$.

Note that when the vector field has a periodic tangency we may still obtain horizontal strips in $In(\vv)$ that are mapped into vertical strips across themselves by the first return map, but there is no guarantee that their image will cross the other strips, and and even less that they will  cross $W^u_{loc}(\ww)$.
\medbreak

\textbf{Case III:} If $P \in  \mathcal{D}$, let $\varphi_0<\varphi_1\in\left(-\pi/2,\pi/2\right)$ be the two solutions of $A(\varphi)-\alpha_\vv E_\ww/\alpha_\ww=0$ where $A(\varphi)$ is the expresssion \eqref{AdePhi} in the proof of Theorem~\ref{prop12} and let
$d=\varphi_1-\varphi_0\in \left(0,\pi/2\right)$.
Then all the positive solutions $\varphi_n$  of $A(\varphi)-\alpha_\vv E_\ww/\alpha_\ww=0$ are of the form
$\varphi_{2n}=\varphi_0+n\pi$, $\varphi_{2n+1}=\varphi_1+n\pi$.
The values of $s\in\left(0,1\right)$ where $dx_\ww/ds=0$ form the decreasing sequence
$s_n=e^{-\frac{\varphi_n}{g_\vv}}e^{\frac{c_2}{g_\vv}}$ that satisfies
$s_n=s_{n-2}e^{-\frac{\pi}{g_\vv}}$ for $n=2,3,\ldots$,  with $\lim_{n\to\infty}s_n=0$.
Since these are the only solutions of $dx_\ww/ds=0$, then $x_\ww(s)$ is monotonic in each interval
$\left(s_{n+1},s_n\right)$.

From  \eqref{eqdxw} it follows that $dx_\ww/ds>0$ if and only if $A(\varphi(s))>C_\vv/\gamma$.
At $\varphi=0$ we have $A(0)=C_\vv a^2$ and this is larger than $C_\vv/\gamma$ because $a>1$ and $\gamma>1$.
Hence $dx_\ww/ds>0$ if $s\in\left(s_1,s_0\right)$.
Since $A(\varphi)$ has period $\pi$, then $dx_\ww/ds>0$ for $s\in\left(s_{2n+1},s_{2n}\right)$.

From Lemma~\ref{lemaxwsn} it follows that $x_\ww(s_{2n+1})-x_\ww(s_{2n})=d$.
As noted at the end of the proof of Theorem~\ref{prop12}, the values of $x_\ww(s_{2n+1})\pmod{2\pi}$ correspond to an irrational rotation around the circle, so these values are dense and  uniformly distributed in $\left[0,2\pi\right)$.
Therefore, given $\tau<d/2$ with $0<\tau<\tau_0$, there exist $n_j\to \infty$ such that
$x_\ww(s_{2n_j+1})\in \left(-\tau,0\right)\pmod{2\pi}$ and hence $x_\ww(s_{2n_j})\in \left(0,\pi/2\right)\pmod{2\pi}$.
Hence, for $s$ in the intervals $\left(s_{2n_j+1},s_{2n_j}\right)$ the curve $x_\ww(s)$ goes across
$\left[-\tau,0\right]$, as required.
\medbreak

\begin{figure}
\begin{center}
\includegraphics[height=4cm]{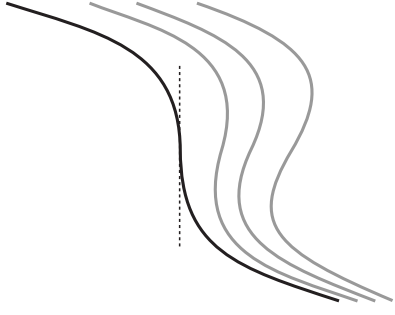}
\end{center}
\caption{\small When the parameter $P$ lies  in $\partial \mathcal{B}$ (\textbf{Case IV} of the proof of Proposition~\ref{propStrips}) the black curve $\eta(\beta(s))$ has a vertical tangent, but $x_\ww(s)$ is still monotonic. However, grey nearby curves
$\eta(\beta_t(s))$, for small $t>0$, may reverse the direction of their turning around $Out(\ww)$.}
\label{figuraInflexao}
\end{figure}

\textbf{Case IV:} If $P \in \partial \mathcal{B}$, then the points where $\frac{dx_\ww}{ds}=0$ are inflection points of $x_\ww(s)$.
Then  the curve $\eta(\beta(s))$ does not reverse the direction of its turning around $Out(\ww)$ and hence $x_\ww$ is a monotonic function of $s$, increasing with $s$ if $\gamma>1$.
However, nearby curves $\eta(\beta_t(s))$ for $t>0$ may reverse their turning at pairs of points near these inflections (see Figure~\ref{figuraInflexao}).
We can choose  $\tau<\tau_0$ and  adapt the arguments of Cases II and III, as appropriate, to obtain inflection points of $x_\ww(s)$ sufficiently far from the interval
$\left[-\tau,0\right]$,  to ensure that the first coordinates of the pairs of turning points
$(\mathrm{mod}\ {2\pi})$  do not fall in that interval.
\Qed

\textbf{Proof of Theorem~\ref{TeoremaTransversal}:}
From Proposition~\ref{propStrips} it follows that  there is a subset of $In(\vv)$ where the first return map is semi-conjugated to a shift in an infinite set of symbols.
Trajectories of points in this set return to $In(\vv)$ infinitely many times, as in assertion~\ref{item6}.

For assertion~\ref{item4}. of Theorem~\ref{TeoremaTransversal} we use the fact that each   horizontal strip of Proposition~\ref{propStrips} is mapped by the first return map $\Psi_{\ww,\vv}\circ \eta$  into a vertical strip across $[0,\tau]\times[0,\tau]$, and hence its image crosses $W^s_{loc}(\vv)$.
Reversing the argument, we get a line of points that come from $W^u_{loc}(\ww)$, intersecting the two lines we get the 2-pulse connections. The $n$-pulse connections may be found iterating this argument or, alternatively, adapting the arguments of \cite{LR, Rodrigues3}.

For assertion~\ref{item5}. we need to show that the first return map is hyperbolic, with a contracting direction.
This is done in Lemma~\ref{LHiperb} below.

Using the restrictions of Proposition~\ref{propStrips} and of Lemma~\ref{LHiperb},
the set $\mathcal{E}$ consists of  those vector fields on $\EU^3$ satisfying (P\ref{H1})--(P\ref{H4}),  not having a periodic tangency, for which $W^s_{loc}(\vv)\cap Out(\ww)$ is close to a vertical line,
with the restrictions $\displaystyle\gamma=\frac{\alpha_\ww}{\alpha_\vv}\frac{C_\vv}{E_\ww} \ne 1$, and
$\displaystyle\delta=\delta_\vv\delta_\ww=\frac{C_\vv C_\ww}{E_\vv E_\ww}>1$.
Clearly, this is an open set  and $\mathcal{E}\cap\mathcal{C}\ne\emptyset$ where $\mathcal{C}$ is  the open set  of Theorem~\ref{teoremaReverteRotacao}.
\Qed

\begin{lemma}\label{LHiperb}
Let $g=\Psi_{\ww,\vv}\circ\eta:In(\vv)\rightarrow In(\vv)$, be the first return map of a vector field on $\EU^3$ satisfying (P\ref{H1})--(P\ref{H4})
and suppose that $\displaystyle\delta=\frac{C_\vv C_\ww}{E_\vv E_\ww}>1$.
Then,
if one of the conditions below hold:
\begin{enumerate}
\item
$P \notin \mathcal{B}$;
\item
$P \in \mathcal{B}$ and $(x,y)$ lies in one of the horizontal strips of Proposition~\ref{propStrips};
\end{enumerate}
 $g$ is hyperbolic at $(x,y)\in In(\vv)$ with $y>0$ sufficiently small.
\end{lemma}
\textbf{Proof:}
From the expressions in Section~\ref{SectionLocalDynamics} we get
$$
\det  Dg(x,y)=c_1^{\delta_\ww} \delta y^{\delta-1} C(\varphi)^{-1+\delta_\ww/2}
\left(1+(c_4-1) g_\ww(a^2-\frac{1}{a^2}) \sin\varphi\cos\varphi \right)
$$
and since $1/a^2\le C(\varphi)\le a^2$,   then $C(\varphi)$ is limited and we  have $\lim_{y\to 0^+} \det Dg(x,y)=0$ if $\delta>1$.
So, for small $y>0$, the derivative   $Dg(x,y)$ has at least one contracting eigenvalue.

The trace of $ Dg(x,y)$ is
$$
 \mathrm{Tr\ }Dg(x,y)=
-c_1^{\delta_\ww} \delta_\ww y^\delta C(\varphi)^{-1+\delta_\ww/2}\left(a^2-\frac{1}{a^2}\right) \sin\varphi\cos\varphi +
\frac{1}{y}\frac{\alpha_\ww}{E_\ww E_\vv C(\varphi)}\left( A(\varphi)-\frac{\alpha_\vv E_\ww}{\alpha_\ww}
\right)
$$
and we want to compute $\lim_{y\to 0} \mathrm{Tr\ }Dg(x,y)$.
Since $\delta>1$, the first summand tends to $0$ as $y$ tends to $0$, and the second summand dominates the limit.
For the second summand, we need to look at the parameters $P$ in the equation.
If $P \notin \mathcal{B}$ then $A(\varphi)-\frac{\alpha_\vv E_\ww}{\alpha_\ww}\ne 0$ has constant sign, and therefore $\lim_{y\to 0} \mathrm{Tr\ }Dg(x,y)=\pm\infty$.
In this case it follows that, for small $y>0$ and any $x$, the derivative  $Dg(x,y)$ has one contracting and one expanding eigenvalue.

For $P \in \mathcal{B}$, the expression  $A(\varphi)-\frac{\alpha_\vv E_\ww}{\alpha_\ww}\ne 0$ does not have constant sign in general, but it does have the same sign inside each one of the horizontal strips of Proposition~\ref{propStrips}.
Without loss of generality, suppose it is positive.
If $(x_n,y_n)$ is any sequence contained in the union of those strips and satisfying
$\lim_{n\to\infty} y_n=0$, then we  have
$\lim_{n\to\infty} \mathrm{Tr\ }Dg(x_n,y_n)=\infty$, and hence  for $(x,y)$ inside a strip
with  small $y>0$, the derivative  $Dg(x,y)$ has one contracting and one expanding eigenvalue.
\Qed

When $P \in \mathcal{B}$, the conclusion of Lemma~\ref{LHiperb} only holds inside the horizontal strips, because they exclude the tangencies of $W^s_{loc}(\vv)$ and $W^u_{loc}(\ww)$.
Near the points where $\beta(s)$ has a vertical tangent, one may find sequences $(x_n,y_n)$ with
$\lim_{n\to\infty} y_n=0$ for which $\lim_{n\to\infty} \mathrm{Tr\ }Dg(x_n,y_n)$ takes any value between $+\infty$ and $-\infty$.
In particular, there are sequences for which $\lim_{n\to\infty} \mathrm{Tr\ }Dg(x_n,y_n)=0$.
At points in such a sequence $Dg(x_n,y_n)$ has two contracting eigenvalues, and this may be additional evidence for the existence of sinks, predicted by Newhouse's results.
Any neighbourhood of the tangency will contain points where the first return map is not hyperbolic, at which one of the eigenvalues crosses the unit circle.

 {The set $\mathcal{D}$ has full Lesbesgue measure.
However, the existence of tangencies of  Theorem \ref{teoremaFinal} does not hold for a full Lesbesgue measure set, only for a  dense subset of vector fields in $\mathcal{D}$.
This is in agreement with Kaloshin's theorem \cite{Kaloshin} on the prevalence of Kupka-Smale systems.}

\section{Example}
\label{numerics}

 In this section we construct a vector field in $\EU^3$ that satisfies properties  (P\ref{H1}) and (P\ref{H4}) and has a connection of one-dimensional invariant manifolds as in (P\ref{H2}).
 As far as we know, no explicit examples of differential equations satisfying (P\ref{H1})--(P\ref{H4}) have been described in the literature, although these conditions follow from the set of properties described by Turaev and Shilnikov \cite{TS98}.
 We present some evidence that the two-dimensional invariant manifolds  intersect transversely, as in (P\ref{H3}), and we use the vector field to obtain numerical simulations that illustrate our results.
Our construction is based on properties of differential equations with symmetry, we
 refer the reader to Golubitsky \emph{et al} \cite{Golubitsky} for more information on the subject.

\subsection{Construction of the example}
We use a
technique presented in Aguiar \emph{et al} \cite{ACL BIF CHAOS} that consists essentially in three steps.
Start with a symmetric vector field on $\RR^3$
with an attracting flow-invariant two-sphere containing a heteroclinic network.
The heteroclinic network involves equilibria and one-dimensional heteroclinic connections that correspond to the intersection of fixed-point subspaces with the invariant sphere.
If one of the symmetries of the vector field in $\RR^3$ is a reflection, then it can be lifted by a rotation to an $\So$--equivariant vector field in  $\RR^4$.
The sphere $\EU^3$ is  flow-invariant and attracting for the lifted vector field,
and a two-sphere of heteroclinic connections arises from
one-dimensional heteroclinic connections  lying outside the plane fixed by the reflection.
Perturbing  the vector field in a way that destroys the $\So$--equivariance and maintains the invariance of the three-sphere breaks the two-dimensional heteroclinic connection into a transverse intersection of invariant manifolds.

The first step in the construction of  \cite{ACL BIF CHAOS,LR3} is to obtain the differential equation in $\RR^3$
\begin{equation}\label{dim3}
\left\{\begin{array}{l}
\dot x=x(1-r^2)-\alpha_1x z+\alpha_2 x z^2\\
\dot y=y(1-r^2)+\alpha_1y z+\alpha_2 y z^2\\
\dot z=z(1-r^2)+\alpha_1(z^2-x^2)-\alpha_2z(x^2+y^2)
\end{array}\right.
\qquad r^2=x^2+y^2+z^2
\end{equation}
that has symmetries
$$
\kappa_1(x,y,z)=(-x,y,z)
\qquad\mbox{and}\qquad
\kappa_2(x,y,z)=(x,-y,z) .
$$
The unit sphere $\EU^2$ is  flow-invariant and globally  attracting, and  $(0,0,\pm 1)$ are equilibria.
From the symmetry it follows that the planes $x=0$ and $y=0$ are flow-invariant, and hence they meet $\EU^2$ in two flow-invariant circles connecting the equilibria $(0,0,\pm 1)$. If  $\alpha_2<0<\alpha_1$
with $\alpha_1+\alpha_2>0$, then these two equilibria are saddles, and there are heteroclinic trajectories going from each equilibrium to the other one, see
 \cite{ACL BIF CHAOS,LR3, Rodrigues}.

Now we adapt the second step in  \cite{ACL BIF CHAOS,LR3, Rodrigues} to obtain property (P\ref{H4}).
Add to \eqref{dim3} a fourth coordinate $\theta$ and the equation $\dot\theta=z$.
Taking $(x,\theta)$ as polar coordinates and rewriting the result in rectangular coordinates
$X=(x_1,x_2,x_3,x_4)=(x\cos\theta,x\sin\theta,y,z)$,
yields a   differential equation in  $\RR^4$ that has the rotational symmetries
$$
(x_1,x_2,x_3,x_4)\to(x_1\cos\varphi-x_2\sin\varphi,x_1\sin\varphi+x_2\cos\varphi,x_3,x_4),
$$
a representation of $\So$.
The unit sphere $\EU^3$ is flow-invariant under the new equation and  attracts every trajectory with non-zero initial condition.
Let $f_0$ be the vector field defined in $\EU^3$ by the new equations.
There are two equilibria given by:
$$
\vv =(0,0,0,+1) \quad \text{and} \quad \ww =(0,0,0,-1)
$$
that, under the conditions on $\alpha_1,\alpha_2$ above, are saddle-foci of different Morse indices.
They share a two-dimensional invariant manifold, $W^s(\vv)=W^u(\ww)$, the two-sphere
$\EU^3\cap\{x_3=0\}$  that lies in the fixed-point subspace
of the symmetry $\tilde\kappa_2(x_1,x_2,x_3,x_4)=(x_1,x_2,-x_3,x_4)$ inherited from $\kappa_2$.
There is also a pair of one-dimensional connections $[\vv\rightarrow\ww]$ in the plane of points fixed by the rotational symmetry.
Since $\dot \theta$ is positive near $\vv$ and negative near $\ww$, the two saddle-foci have different chirality. The vector field $f_0$  satisfies properties  (P\ref{H1}), (P\ref{H2}) and (P\ref{H4}).

The third and final step is to obtain property  (P\ref{H3}).
This implies  breaking the rotational symmetry.
The vector fields $f_\lambda=f_0+\lambda g$ do not have the two-dimensional connection
for generic $g$ and for small values of $\lambda$.
Let $g$ be a vector field in $\RR^4$  tangent to $\EU^3$,
that does not have the symmetry $\tilde{\kappa}$ nor all the rotational symmetries $\So$,
but for which $\tilde\kappa_1(x_1,x_2,x_3,x_4)=(-x_1,-x_2,x_3,x_4)$ is still a symmetry.
The last requirement ensures that the  one-dimensional connection remains for the perturbed vector field $f_\lambda=f_0+\lambda g$,  maintaining the other properties.
An example of the result of this construction  is equation \eqref{example} below.

Our results are illustrated with numerical simulations, which have been obtained using the dynamical systems package \emph{Dstool}.

\subsection{The  example}

Our example is the one-parameter family $f_\lambda(X)$ of vector fields on $\EU^3\subset\RR^{4}$, defined
by  the differential equation in $\RR^{4}$:
\begin{equation}
\left\{
\begin{array}{l}
\dot{x}_{1}=x_{1}(1-r^2)-x_4 x_2-\alpha_1x_1x_4+\alpha_2x_1x_4^2\\
\dot{x}_{2}=x_{2}(1-r^2)+x_4 x_1-\alpha_1x_2x_4+\alpha_2x_2x_4^2 \\
\dot{x}_{3}=x_{3}(1-r^2)+\alpha_1x_3x_4+\alpha_2x_3x_4^2+\lambda x_1x_2x_4\\
\dot{x}_{4}=x_{4}(1-r^2)-\alpha_1(x_3^2-x_1^2-x_2^2)-\alpha_2x_4(x_1^2+x_2^2+x_3^2)-\lambda x_1x_2x_3 \\
\end{array}
\label{example}
\right.
\end{equation}
where  $r^2=x_{1}^{2}+x_{2}^{2}+x_{3}^{2}+x_{4}^{2}$ and $\alpha_2<0<\alpha_1$ with $\alpha_1+\alpha_2>0$.
\medbreak
The unit sphere $\EU^3$ is  invariant under the flow of \eqref{example} and every trajectory with nonzero initial condition is asymptotic to it in forward time, so $f_\lambda(X)$ is a well defined  vector field  on $\EU^3$  for each $\lambda$, with the two equilibria $\vv$ and $\ww$.
The linearisation of $f_\lambda(X)$
at $(0,0,0,\varepsilon)$ with
$\varepsilon=\pm 1$ has non-radial eigenvalues  $$\alpha_2-\varepsilon\alpha_1\pm i \quad \text{and} \quad
\alpha_2+\varepsilon\alpha_1.$$
Under the conditions above, $\vv$ and $\ww$ are hyperbolic saddle-foci,
$\vv$ has one-dimensional unstable  manifold and two-dimensional stable manifold; $\ww$ has one-dimensional stable manifold and two-dimensional unstable manifold.

\begin{figure}
\begin{center}
\includegraphics[height=7cm]{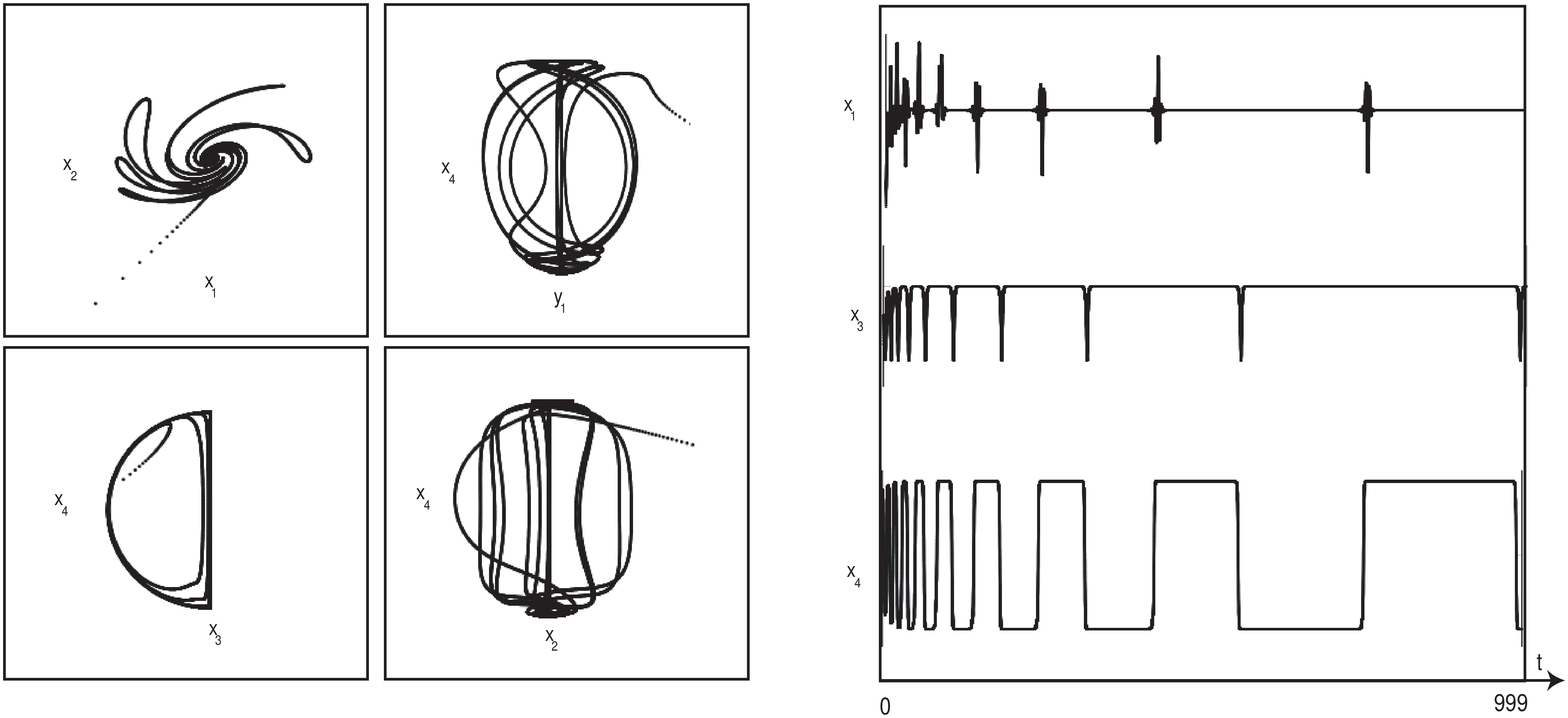}
\end{center}
\caption{\small Example of a solution of equation \eqref{example} for $\lambda=0$ that accumulates asymptotically on one of the Bykov cycles in the network. Left: Projection in the $(x_1, x_2)$, $(x_1, x_4)$, $(x_3, x_4)$ and $(x_2, x_4)$--planes of the trajectory with initial condition $(-0.5000, -0.1390, -0.8807, 0.3013)$, corresponding to $\alpha_1=1$ and $\alpha_2=-0.1$. Right: Time series for the same trajectory.}
\label{unperturbed}
\end{figure}

For $\lambda=0$, the one-dimensional invariant manifolds of $\vv$ and $\ww$ lie in the invariant circle $Fix(\textbf{SO(2)}) \cap \EU^3$ and the two-dimensional  invariant manifolds  lie in the invariant two-sphere $Fix(\ZZ_2(  \tilde{\kappa}_2 ))\cap \EU^3$.
Thus, symmetry forces the invariant manifolds of $\vv$ and $\ww$ to be in a very special position: they coincide.
The two saddle-foci, together with their invariant manifolds form a heteroclinic network
$\Sigma$ that is asymptotically stable by the Krupa and Melbourne  criterion \cite{Krupa e Melbourne 1, Krupa e Melbourne}.
Indeed since $\alpha_2<0<\alpha_1$, it follows that
$
\delta=\frac{C_\vv}{E_\vv}\frac{C_\ww}{E_\ww}=\left(\frac{\alpha_2-\alpha_1}{\alpha_2+\alpha_1}\right)^2>1,
$
where $E_X$ and $C_X$ denote the real parts of the expanding and contracting eigenvalues of
$Df_0(X)$ at $X=\vv$ and $X=\ww$, respectively. The network $\Sigma$ can be decomposed into two cycles. Due to the symmetry and to the asym\-ptotic stability, trajectories whose initial condition lies outside the invariant fixed point subspaces will approach in positive time one of the cycles. The fixed point hyperplane $Fix(\ZZ_2(  \tilde{\kappa}_2 ))$ prevents random visits to the two cycles; a trajectory that approaches one of the cycles in $\Sigma$. The time-series of  Figure~\ref{unperturbed} shows the increasing intervals of time spent near the equilibria. The sojourn time in any neighbourhood of one of the saddle-foci increases geometrically with ratio $\delta$.

\medbreak

\begin{figure}
\begin{center}
\includegraphics[height=7cm]{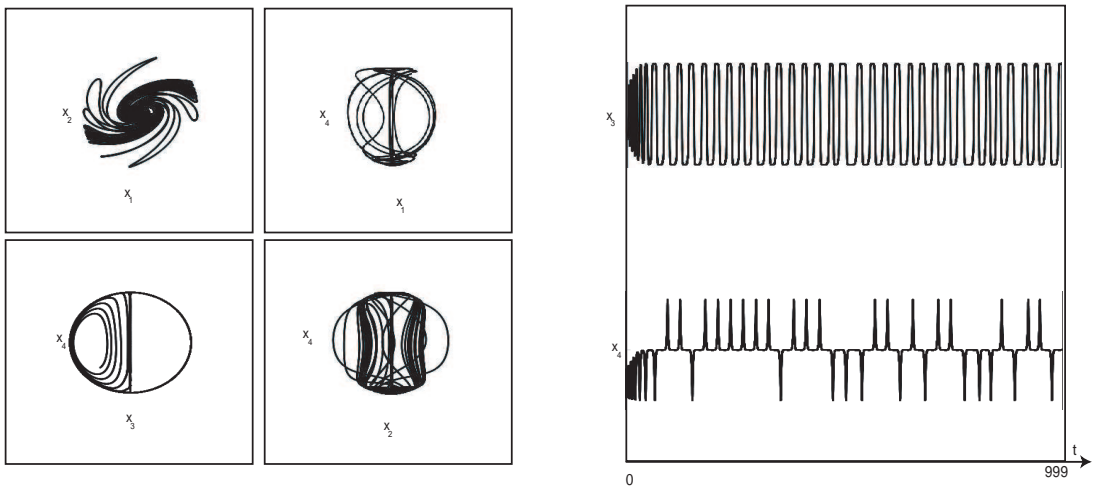}
\end{center}
\caption{\small Example of a solution of  equation \eqref{example} for $\lambda=0.05$ that visits the two primary Bykov cycles in the network. If the two-dimensional manifolds $W^s(\vv)$ and $W^u(\ww)$ first meet transversely (which is generically the case), then there are trajectories doing these visits in any prescribed order in a behaviour called switching. Left: Projection in the $(x_1, x_2)$, $(x_1, x_4)$, $(x_3, x_4)$ and $(x_2, x_4)$--planes of the trajectory with initial condition $(-0.5000, -0.1390, -0.8807, 0.3013)$, corresponding to $\alpha_1=1$ and $\alpha_2=-0.1$. Right: Time series for the same trajectory.}

\label{perturbed}
\end{figure}

When $\lambda \neq 0$, according to the terminology of \cite{Robert}, an explosion on the non-wandering set occurs. The parameter $\lambda$ should control the transversality of the 2-dimensional local invariant  manifolds. Care needs to be taken with numerical integration of systems with heteroclinic
cycles and networks, because rounding errors may
cause qualitatively incorrect results. We have not attempted to prove analytically the transversality -- we defer this analysis to a future paper.
The switching  mechanism described in \cite{ACL NONLINEARITY, ALR}
operating in our network ensures that most trajectories will visit
most parts of a neighbourhood of the network, as suggested in Figure~\ref{perturbed}.

\subsection{Different chirality}
The example given here in (\ref{example}) is similar to that reported in  \cite{LR3}.
The latter has been constructed using the standard lift technique and thus the equation for  the angular coordinate $\theta$ in the plane $(x_1,x_2,0,0)$ is
$\dot{\theta}=1$. Since this plane is perpendicular to the plane where the connection  $[\vv \rightarrow \ww]$ lies, trajectories  must turn around the connection in the same direction and the nodes have the same chirality. In the case of example (\ref{example}), the chirality at the two saddle-foci is different, also by construction: since $\dot{\theta}=z$,
near $\vv$  we have $\dot\theta>0$ and near $\ww$ we have $\dot\theta<0$.
Chirality does not have any impact on the asymptotic stability of the network, making our work completely consistent with the Krupa and Melbourne criterion \cite{Krupa e Melbourne 1,Krupa e Melbourne}.
In the simulations presented here we have used the same parameters and initial conditions reported in \cite{LR3}, to facilitate the comparison of the two cases.

\section{Chirality in Bykov's results}
We conclude with a brief discussion to explain that 
Bykov implicitly assumed that the chiralities of the two nodes are different in formulas (3.1),(3.4),(3.5) of \cite{Bykov}.
Instead of looking at the walls of the cylinders (as we do in Theorem \ref{prop12}), he studies the intersection of the two-dimensional invariant manifolds of the saddle-foci with the disks that we call $Out(\vv)$ and $In(\ww)$, described by his formula (3.1).
The two-dimensional invariant manifolds intersect these disks on spirals.
After this, Bykov compares the spiral  that we would denote $W^s_{loc}(\vv)\cap In(\ww)$ to the deformed spiral image of the spiral $W^u_{loc}(\ww)\cap Out(\vv)$ by  a transition map similar to our $\Psi_{\vv,\ww}$.
In his notation the non-real eigenvalues at the saddle-foci are $\alpha_j\pm i\omega_j$, $j=1,2$ and the angular coordinates in $Out(\vv)$ and $In(\ww)$ are denoted $\varphi_j$.
The formulas that represent the evolution of the angular component  $\xi_j$ of the spirals are:
$$
\xi_1=d_1 e^{-\varphi_1/\omega_1}(1+\chi_{21}(0, \varphi_1,0))
\quad
\text{and}
\quad
\xi_2=d_2 e^{-\varphi_2/\omega_2}(1+\chi_{22}(0, \varphi_2,0))
$$
for some constants $d_j$ and some maps $\chi_{2j}$.
 Indeed, observe the same sign of the exponent of $e$ in the formulas above, that corresponds to different chirality, since one formula is obtained by looking directly at the flow and the other corresponds to the flow taken in reverse time.
 The same happens in Formula (3.1) of \cite{Bykov99}. 
 
 Bykov never comments on the chiralities of the nodes, assuming implicitly that they are different. 
 Therefore, in the cross sections, the spirals corresponding to the two-dimensional invariant manifolds of the saddle-foci are oriented in the same way. 
 This explains the orientation of the spirals of Figure 2 of \cite{Bykov} in contrast to those shown in Figure 11 of \cite{KLW}, that turn in opposite directions because the nodes have the same chirality.
Bykov's condition $G<0$ in Theorem 3.2 of \cite{Bykov} is analogous to our condition \eqref{final_formulae} that defines the set
$\mathcal{B}$ where tangencies are dense.


\begin{thebibliography}{99}


\bibitem{ACL NONLINEARITY} M.A.D. Aguiar, S.B.S.D. Castro, I. S. Labouriau,
\emph{Dynamics near a heteroclinic network,} Nonlinearity, 18, 391--414, 2005

\bibitem{ACL BIF CHAOS} M.A.D. Aguiar, S.B.S.D. Castro, I. S. Labouriau,
\emph{Simple Vector Fields with Complex Behaviour}, Int. Jour. of
Bifurcation and Chaos, Vol. {16},  2, 369--381,  2006

\bibitem{ALR} M.A.D. Aguiar, I.S. Labouriau, A.A.P. Rodrigues, \emph{Switching near a heteroclinic network of rotating nodes}, Dynamical Systems: an International Journal, Vol. 25(1), 75--95, 2010


\bibitem{AC} P. Ashwin, P. Chossat, \emph{Attractors for Robust Heteroclinic Cycles with Continua of Connections}, J. Nonlinear Sci., Vol. {8}, 103--129, 1998

\bibitem{Bowen} R. Bowen, \emph{A horseshoe with positive measure}, Invent. Math. 29, 203--204, 1975

\bibitem{Bowen75} R. Bowen, \emph{Equilibrium States and the Ergodic Theory of Anosov Diffeomorphisms}, Lect. Notes in Math, 1975

\bibitem{Bykov99}  V. V. Bykov, \emph{On systems with separatrix contour containing two saddle-foci}, Journal of Mathematical Sciences 95, 2513--2522, 1999

\bibitem{Bykov} V. V. Bykov, \emph{Orbit Structure in a Neighbourhood of a Separatrix Cycle Containing Two Saddle-Foci}, Amer. Math. Soc. Transl, Vol. 200, 87--97, 2000


\bibitem{Delshams} A. Delshams, S. Gonchenko, M. Gonchenko, J. L\'azaro, \emph{Mixed dynamics of two-dimensional reversible maps with a symmetric couple of quadratic homoclinic tangencies},  arXiv:1412.1128, 2014

\bibitem{DumortierEtAl} F. Dumortier, S. Iba\~nez, H. Kokubu, \emph{Cocoon bifurcation in three-dimensional reversible vector fields}, Nonlinearity 19, 305--328, 2006

\bibitem{Field book} M. Field, \emph{Lectures on bifurcations, dynamics and symmetry}, Pitman Research Notes in Mathematics Series, Vol. {356}, Longman, 1996

\bibitem{GS} N.K. Gavrilov, L.P. Shilnikov, \emph{On three-dimensional dynamical systems close to systems with a structurally unstable homoclinic curve}, Part I, Math. USSR, Sbornik 17, 467--485, 1972


\bibitem{T-points Glendinning} P. Glendinning, C. Sparrow, \emph{T-points: A codimension Two Heteroclinic Bifurcation}, J. Stat. Phys., 43,  No. 3--4, 479--488, 1986

\bibitem{Golubitsky} M.I. Golubitsky, I. Stewart,  D.G. Schaeffer, \emph{Singularities and Groups in Bifurcation Theory }, Vol. {II}, Springer, 2000

\bibitem{Gonchenko} S.V. Gonchenko,  M.C. Li, \emph{On hyperbolic dynamics of multidimensional systems with homoclinic tangencies of arbitrary orders}, preprint, 2011

\bibitem{GSTchaos} S. V. Gonchenko, L.P. Shilnikov, D.V. Turaev,
\emph{Dynamical phenomena in systems with structurally unstable Poincare homoclinic orbit},
Chaos 6(1), 15--31, 1996

\bibitem{GST} S.V. Gonchenko, L.P. Shilnikov, D.V. Turaev, \emph{Homoclinic tangencies of arbitrarily high orders in conservative and dissipative two-dimensional maps}, Nonlinearity 20, 241--275, 2007

\bibitem{Guckenheimer and Holmes} J. Guckenheimer, P. Holmes, \emph{Nonlinear Oscillations, Dynamical Systems and Bifurcations of Vector Fields}, Springer-Verlag, 1983


\bibitem{Hayashi} S. Hayashi, \emph{Connecting invariant manifolds and the solution of the $C^1$-stability and $\Omega$-stability conjectures for flows}, Ann. Maths, 145, 81Ð137, 1997



\bibitem{HS} A.J. Homburg, B. Sandstede, \emph{Homoclinic and Heteroclinic Bifurcations in Vector Fields}, Handbook of Dynamical Systems, Vol. 3, North Holland, Amsterdam, 379--524, 2010

\bibitem{Kaloshin} V. Kaloshin, \emph{Some prevalent properties of smooth dynamical systems}, 
Trudy Matematichskogo Instituta imeni V.A. Steklova, 213, 123--151, 1996, Translation in Proc. Steklov Inst. Math., 213, 115--140, 1996


\bibitem{Katok2} A. Katok, \emph{Lyapunov exponents, entropy and periodic orbits for diffeomorphisms}, IHES Publ. Math.,  Vol. 51, 137--173, 1980

\bibitem{KLW} J. Knobloch, J.S.W. Lamb, K. N. Webster, \emph{Using Lin's method to solve Bykov's problems}, J. Diff. Eqs., 257(8), 2984--3047, 2014

%


\bibitem{Krupa e Melbourne 1} M. Krupa, I. Melbourne, \emph{Asymptotic Stability of Heteroclinic Cycles in Systems with Symmetry, }Ergodic Theory
and Dynam. Sys., Vol. {15}, 121--147, 1995

\bibitem{Krupa e Melbourne} M. Krupa, I. Melbourne, \emph{Asymptotic Stability of Heteroclinic Cycles in Systems with Symmetry,\ II, }Proc. Roy.
Soc. Edinburgh, 134A, {1177--1197},  2004


\bibitem{LR} I.S. Labouriau, A.A.P. Rodrigues, \emph{Global Generic Dynamics Close to Symmetry},  Journal of Differential Equations, Vol. 253, 8  2527--2557, 2012

\bibitem{LR2}  I.S. Labouriau, A.A.P. Rodrigues, \emph{Partial Symmetry Breaking and Heteroclinic Tangencies}, in S. Ib\'a\~nez,
J.S. P\'erez del R\'io, A. Pumari\~no and J.A. Rodr\'iguez (Eds), Progress and challenges in dynamical systems,
Proceedings in Mathematics and Statistics; Springer-Verlag, 281--299, 2013

\bibitem{LR2015} I.S. Labouriau, A.A.P. Rodrigues, \emph{Global Bifurcations Close to Symmetry}, arXiv: 1504.01659v1


\bibitem{LS} J. Lamb, O. Stenkin, \emph{Newhouse regions for reversible systems with infinitely many stable, unstable and elliptic periodic orbits}, Nonlinearity, 17, 1217--1244, 2004

\bibitem{LTW} J. Lamb, M. Teixeira, K. Webster, \emph{Heteroclinic bifurcations near Hopf-zero bifurcation in reversible vector fields
in $\RR^3$}, J. Differential Equations 219, 78--115, 2005

%


\bibitem{MV} L. Mora, M. Viana, \emph{Abundance of strange attractors}, Acta Math. 171, 1--71, 1993

\bibitem{Newhouse2} S.E. Newhouse, \emph{Diffeomorphisms with infinitely many sinks}, Topology, {13}, 9--18, 1974

\bibitem{Newhouse1} S.E. Newhouse, \emph{The abundance of Wild Hyperbolic Sets and Non-Smooth Stable Sets for Diffeomorphisms}, Publ. Math. Inst. Hautes \'Etudes Sci., Vol. {50},  101--151, 1979

\bibitem{OS} I.M. Ovsyannikov, L.P. Shilnikov, \emph{On systems with saddle-focus homoclinic curve}, Math. USSR
Sbornik, 58,  557--574, 1987


\bibitem{PT} J. Palis, F. Takens, \emph{Hyperbolicity and sensitive chaotic dynamics at homoclinic bifurcations}, Cambridge University Press, 1993

\bibitem{Robert} C. Robert, K. Alligood, E. Ott, J. Yorke, \emph{Explosions of chaotic sets}, Physica D, 44--61, 2000


\bibitem{Rodrigues3} A.A.P. Rodrigues, \emph{Repelling dynamics near a Bykov cycle}, Journal of Dynamics and Differential Equations, Vol. 25 (3), 605--625, 2013


\bibitem{LR3} A.A.P. Rodrigues, I.S. Labouriau, \emph{Spiralling dynamics near  heteroclinic networks},  Physica D: Nonlinear Phenomena, 268, 34--49, 2014

\bibitem{Rodrigues} A.A.P. Rodrigues, I.S. Labouriau, M.A.D. Aguiar, \emph{Chaotic Double Cycling}, Dynamical Systems: an International Journal, Vol. 26(2), 199--233, 2011

\bibitem{Samovol} V.S. Samovol, \emph{Linearization of a system of differential equations in the neighbourhood of a singular point, }Sov.\
Math. Dokl, Vol. {13}, 1255--1959, 1972



\bibitem{Shilnikov_67} L.P. Shilnikov, \emph{On a Poincar\'e--Birkhoff problem}, Math. USSR Sb. 74(3), 353--371, 1967



\bibitem{TS98} D. Turaev, L.P. Shilnikov, \emph{An example of a wild strange attractor}, Math. USSR Sb. 189 (2), 353--371, 1967

\bibitem{Wen} L. Wen, \emph{Homoclinic tangencies and dominated splittings}, Nonlinearity 15, 1445--1469, 2002

\bibitem{YA} J.A. Yorke, K.T. Alligood, \emph{Period-Doubling Cascade of Attractors: A prerequisite for Horseshoes}, Communications in Mathematical Physics, 101, 305--321, 1985

\end{thebibliography}
\end{document}